\documentclass[12pt]{article}

\usepackage{cancel,soul} %
\usepackage[latin1]{inputenc}
\usepackage{bbm}
\usepackage{stmaryrd}
\usepackage{latexsym}
\usepackage{amsfonts}
\usepackage{amsmath,amssymb,amscd}
\usepackage{yfonts}
\usepackage{dsfont}
\usepackage{mathabx}
\usepackage[dvips]{graphicx}
\usepackage{epsfig}
\usepackage{psfrag}
\usepackage[font={small}]{caption}
\usepackage{color}
\usepackage{float}
\usepackage{upgreek}
\usepackage{hyperref}

\usepackage{color,xcolor}
\definecolor{dviolet}{RGB}{102,0,204}

\DeclareFontFamily{U}{txsyc}{}
\DeclareFontShape{U}{txsyc}{m}{n}{
   <-> txsyc%
}{}
\DeclareFontShape{U}{txsyc}{bx}{n}{
   <-> txbsyc%
}{}
\DeclareFontShape{U}{txsyc}{l}{n}{<->ssub * txsyc/m/n}{}
\DeclareFontShape{U}{txsyc}{b}{n}{<->ssub * txsyc/bx/n}{}
\DeclareSymbolFont{symbolsC}{U}{txsyc}{m}{n}
\SetSymbolFont{symbolsC}{bold}{U}{txsyc}{bx}{n}
\DeclareFontSubstitution{U}{txsyc}{m}{n}
\DeclareMathSymbol{\df}{\mathrel}{symbolsC}{"42}
\DeclareMathSymbol{\fd}{\mathrel}{symbolsC}{"43}
\DeclareMathSymbol{\lJoin}{\mathrel}{symbolsC}{"58}
\DeclareMathSymbol{\rJoin}{\mathrel}{symbolsC}{"59}

\newcommand{\cA}{{\cal A}}
\newcommand{\cB}{{\cal B}}

\newcommand{\cE}{{\cal E}}
\newcommand{\cF}{{\cal F}}

\newcommand{\cH}{{\cal H}}

\newcommand{\cM}{{\cal M}}

\newcommand{\EE}{\mathbb{E}}

\newcommand{\LL}{\mathbb{L}}
\newcommand{\NN}{\mathbb{N}}
\newcommand{\PP}{\mathbb{P}}

\newcommand{\RR}{\mathbb{R}}
\renewcommand{\SS}{\mathbb{S}}

\newcommand{\ZZ}{\mathbb{Z}}

\newcommand{\frd}{\mathfrak{d}}

\newcommand{\iy}{\infty}
\newcommand{\lt}{\left}
\newcommand{\me}{\medskip}

\newcommand{\pa}{\partial}
\newcommand{\ri}{\rightarrow}
\newcommand{\rt}{\right}
\newcommand{\sm}{\smallskip}

\newcommand{\wi}{\widetilde}
\newcommand{\wit}{\widehat}
\newcommand{\card}{\mathrm{card}}

\newcommand{\fo}{\forall\ }

\newcommand{\lve}{\lt\vert}
\newcommand{\lVe}{\lt\Vert}

\newcommand{\rve}{\rt\vert}
\newcommand{\rVe}{\rt\Vert}

\newcommand{\St}{\,:\,} 

\newcommand{\un}{\mathds{1}}

\newcommand{\Vect}{\mathrm{Vect}}
\newcommand{\vvv}{\vert\!\vert\!\vert}
\newcommand{\bq}{\begin{eqnarray*}}
\newcommand{\bqn}[1]{\begin{eqnarray}\label{#1}}
\newcommand{\eq}{\end{eqnarray*}}
\newcommand{\eqn}{\end{eqnarray}}
\newcommand{\wwtbp}{\par\hfill $\blacksquare$\par\me\noindent}
\newcommand{\thistitlepagestyle}{}
\newcommand{\lin}{\llbracket}
\newcommand{\rin}{\rrbracket}

\newcommand{\ttsim}{\raise.17ex\hbox{$\scriptstyle\mathtt{\sim}$}}

\newtheorem{pro}{Proposition} 

\newtheorem{lem}[pro]{Lemma}
\newtheorem{theo}[pro]{Theorem}

\renewcommand{\thepro}{\arabic{pro}}

\newenvironment{rem}
{\par\me\refstepcounter{pro}\noindent{\bf Remark \thepro\ }}
{\par\hfill $\square$\par\sm\noindent}
\newenvironment{remm}[1]
{\par\me\refstepcounter{pro}\noindent{\bf Remark \thepro\ #1}}
{\par\hfill $\square$\par\sm\noindent}

\newcommand{\proof}{\par\me\noindent\textbf{Proof}\par\sm\noindent}
\newcommand{\prooff}[1]{\par\me\noindent\textbf{#1}\par\sm\noindent}

\newcommand{\cc}{^{\mathrm{c}}}

\newcommand{\rr}{^{(r)}}
\newcommand{\uu}{^{(1)}}
\newcommand{\umu}{\underline{\smash{\mu}}}
\newtheorem{example}{Example}
\newtheorem{theorem}{Theorem}

\newtheorem{prop}{Proposition}

\newcommand{\Int}{{\rm Int}}
\newcommand{\R}{\mathbb{R}}
\newcommand{\N}{\mathbb{N}}

\newcommand{\A}{\mathcal{A}}

\newcommand{\FS}{{\sf FS}}
\newcommand{\MS}{{\sf MS}}
\newcommand{\RM}{{\sf RM}}

\newcommand{\diver}{{\rm div}\,}

\newcommand{\vertiii}[1]{{\left\vert\kern-0.25ex\left\vert\kern-0.25ex\left\vert #1 
    \right\vert\kern-0.25ex\right\vert\kern-0.25ex\right\vert}}

\title{Higher order Cheeger inequalities for Steklov eigenvalues}
\date{}
\author{Asma Hassannezhad${}^\dagger$ and Laurent Miclo${}^\ddagger$
}
\date{\vbox{\copy0
\vskip5mm
\copy1
}
}
\begin{document}

\setbox0=\vbox{
\large
\begin{center}
 ${}^\dagger$Institut Mittag-Leffler, Sweden
\end{center}
}
 \setbox1=\vbox{
 \large
 \begin{center}
 ${}^\ddagger$Institut de Mathématiques de Toulouse, UMR 5219\\
Université de Toulouse and CNRS, France\\
 \end{center}
}
\setbox5=\vbox{
\hbox{${}^\dagger$hasannezhad@mittag-leffler.se\\}
\vskip1mm
\hbox{Institut Mittag-Leffler\\}
\hbox{Aurav\"agen 17\\}
\hbox{SE-182 60 Djursholm, Sweden\\}
}

\setbox6=\vbox{
\hbox{${}^\ddagger$miclo@math.univ-toulouse.fr\\}
\vskip1mm
\hbox{Institut de Mathématiques de Toulouse\\}
\hbox{Université Paul Sabatier\\}
\hbox{118, route de Narbonne\\}
\hbox{31062 Toulouse Cedex 9, France\\}
}

\maketitle
\thistitlepagestyle
\begin{abstract}
We prove a lower bound for the $k$-th Steklov eigenvalues in terms of an isoperimetric constant called the $k$-th Cheeger-Steklov constant in three different situations:  finite spaces,  measurable spaces, and Riemannian manifolds. These lower bounds can be considered as  higher order Cheeger type inequalities for the Steklov eigenvalues. In particular it extends the Cheeger type inequality for the first nonzero Steklov eigenvalue previously studied by Escobar in 1997 and by Jammes in 2015 to higher order Steklov eigenvalues. 
The technique we develop to get this lower bound is based on considering a family of accelerated Markov operators in the finite and mesurable situations and of mass concentration deformations of the Laplace-Beltrami operator
in the manifold setting which converges uniformly to  the Steklov operator.  As  an intermediary step in the proof of the higher order Cheeger type inequality, we define  the Dirichlet--Steklov connectivity spectrum  and show that  the  Dirichlet connectivity spectra of this family of operators converges to (or is bounded by)  the Dirichlet--Steklov spectrum uniformly.   Moreover, we obtain bounds for the Steklov eigenvalues in terms of its Dirichlet-Steklov connectivity spectrum which is interesting in its own right and is more robust than the higher order Cheeger type inequalities.   The Dirichlet--Steklov spectrum  is closely related to the Cheeger--Steklov constants. 
\end{abstract} 
\null
{\small

}\par

\section{Introduction}
Let $(M,g)$ be a compact Riemannian manifold of dimension $n$ with smooth boundary, the Steklov eigenvalue problem is \begin{equation}\label{a-stek}\left\{\begin{array}{ll}
\Delta f=0,&\quad\hbox{in}\;\; M\\
\frac{\partial f}{\partial \nu}=\sigma f,&\quad \hbox{on} \;\; \partial M
\end{array}\right.\end{equation}
where $\Delta=\diver\nabla$ is the Laplace--Beltrami operator on $M$ and $\nu$ is the unit outward normal vector along $\partial M$. Its spectrum consists of a sequence of nonnegative real numbers with  accumulation point only at infinity. 
We denote the sequence of the Steklov eigenvalues by $$0=\sigma_1\le\sigma_2\le\cdots\le\sigma_k\le\cdots\nearrow \infty$$
The Steklov eigenvalues can be also considered as the eigenvalues  of the Dirichlet--to--Neumann operator $$
 S:C^{\infty}(\partial M)\to C^\infty(\partial M)$$
 $$f\mapsto \frac{\partial F}{\partial \nu}$$
 where $F$ is the harmonic extension of $f$ into the interior of $M$. The Steklov problem was  first introduced
by  Steklov  \cite{stek02} in 1902  for bounded domains of the plane. 
Many interesting developments and progress in the study of the Steklov  problem have been attained in recent years. We refer the reader to the survey paper \cite{GP} and the references therein for  recent developments, and to \cite{MR3137253} for a historical account. 
The relationship between the Steklov eigenvalues and  geometry of the underlying space, and also its similarity and difference with the Laplace eigenvalues have been a main focus of interest and a source of inspiration, see for example \cite{MR1479552,MR2770439,MR2807105,MR2838029,MR1696453,MR3461367,MR3449183}.\\ 

The focus of this paper is on obtaining lower bounds for the $k$-th Steklov eigenvalue $\sigma_k$ in terms of some isoperimetric constants in three different settings. Our results can be viewed as counterparts of the higher order Cheeger inequalities   for the Laplace eigenvalues in discrete setting proved by Lee,
Oveis Gharan and 
Trevisan~\cite{MR2961569}, and in  manifold setting by the second author~\cite{MR3323580}.  It is also an extension of Escobar's  \cite{MR1479552,MR1696453} and Jammes' \cite{MR3449183} results for $\sigma_2$. We first recall previous results known in this direction.   \\
Let $\A$ denote the family of all  nonempty open subsets $A$ of $M$ with piecewise smooth boundary. For every $A\in\A$ let $\mu(A)$ denote its  Riemannian measure and $\umu(\partial A)$ denote the $(n-1)$-dimensional Riemannian measure of $\partial A$. We define for every $A\in \A$ the isoperimetric ratios 
\bqn{bdry3}
\eta(A)&\df& \frac{\umu(\pa_i A)}{\mu( A)}\\
\eta'(A)&\df& \frac{\umu(\pa_iA)}{\umu(  \bar A\cap \partial M)}
\eqn
where $\partial_i A:=\partial A\cap {\Int\, M}$. Here $\Int\, M$ denotes the interior of $M$.  Consider the following isoperimetric constants
$$h_2(M):=\inf_{A}\max\{\eta(A),\eta(M\setminus A)\}$$
$$h'_2(M):=\inf_{A}\max\{\eta'(A),\eta'(M\setminus A)\}$$
The constant  $h_2(M)$ is the well-known Cheeger constant \cite{MR0402831}. 
Motivated by the celebrated result of Cheeger \cite{MR0402831}, Escobar  \cite{MR1479552,MR1696453} introduced the isomerimetric constant  $h_2'(M)$ and obtained a lower bound for $\sigma_2$  in terms of this isoperimetric constant  and the  first nonzero eigenvalue of a Robin problem. Recently, Jammes~\cite{MR3449183} obtained a simpler and more explicit lower bound for $\sigma_2$ in terms of  an isoperimetric $\tilde h_2'(M)$ similar to the one introduced by Escobar, and the Cheeger constant $h_2(M)$: 
\bqn{jammes1}\sigma_2(M)&\ge&\frac{1}{4} \tilde{h}'_2(M)h_2(M)\eqn
where $ \tilde{h}'_2(M):=\inf\left\{\eta'(A)\,:\,A\in \A,\,\,\hbox{and}\,\,\mu(A)\le\frac{\mu(M)}{2}\right\}$. The proof of (\ref{jammes1}) is simple and only uses the co-area formula.  The constants $h'_2(M)$ and $\tilde{h}'_2(M)$  are interesting geometric quantities. It is an intriguing question if similar geometric lower bounds hold for higher order Steklov eigenvalues $\sigma_k$. We give an affirmative answer to this question not only in Riemannian setting but also in the setting of finite and measurable spaces. \\
 Let $(M,\mu)$ be a measure space and $V$ a proper subset of $M$, and let $L$ be an operator acting on a functional subspace $\mathcal H$ of $\LL^2(\mu)$. Throughout the paper we deal with either of three different settings listed below: 
 \begin{itemize}
 \item[(\FS)] Finite state spaces: $M$ is a finite set, $V$ is a proper subset of cardinality $v$, $L$ is a reversible irreducible Markov generator  and $\mu$ is its unique invariant probability measure. Here $\mathcal H$ is the space of functions on $M$  denoted by $\cF(M)$.
 \item[(\MS)] Measurable state spaces: $(M,\mu)$ is a probability measure space  with $\sigma-$algebra $\mathcal M$, and  $V$ is a measurable subset of $M$ such that $0<\mu[V]<1$. Here,  $L$ is a Markov generator of the form $P-I$, where $P$ is a Markov kernel
reversible with respect to $\mu$ and $I$ is the identity, and $\mathcal{H}=\LL^2(\mu)$. 
 \item[(\RM)] Riemannian manifolds: $M$ is a compact  Riemannian manifold with smooth boundary $\partial M$, $\mu$ is its Riemannian measure, $L$ is the Laplace-Beltrami operator $\Delta$, and $\mathcal H$ is the Sobolev space $H^1(\mu)$. Here $V$ is equal to  $\partial M$.
 \end{itemize} 
 With the help of $L$ we define an operator  $S$ on $V$ and call it the Steklov operator. In setting (\RM), the operator $S$ we consider is in fact the Dirichlet--to--Neumann operator discussed above. For the definition of $S$ in  (\FS)  and (\MS) settings we  refer to definitions \eqref{defFS}  in Section \ref{FS}, and \eqref{Skern} in Section \ref{MS}, respectively. We denote  the eigenvalues of $S$ by $\sigma_k(M)$ or simply $\sigma_k$.  Let $\A$ be a family of admissible sets in $M$:
 \begin{itemize} \item in (\FS) settings, $\A$ is the set of all nonempty subsets of $M$;
\item in (\MS) setting, $\A$ is the set of all non-negligible elements of $\mathcal M$, i.e.  $A\in \mathcal M$ such that~$0<\mu[A]\le1$;
\item in (\RM) setting, $\A$ is the set of all  nonempty open domains $A$ in $M$ such that $\partial_eA:=\bar A\cap \partial M$ and $\partial_i A:=\partial A\cap M$ are smooth manifolds of  dimension $n-1$ when they are nonempty. \end{itemize}
In (\FS) and (\MS) settings, we introduce
the boundary of any $A\in\cA$ via
\bq
\pa A&\df& \{(x,y)\St x\in A, \, y\in A\cc\}\eq

 and define the following isoperimetric ratios
 \bq
\eta(A)&\df& \frac{\umu(\pa A)}{\mu( A)}\\
\eta'(A)&\df& \frac{\umu(\pa A)}{\mu(  A\cap V)}
\eq
where $\umu$ is a measure on $M\times M$. We refer to \eqref{bdry1} and \eqref{bdry2} for the definition of $\umu$ in  (\FS) and (\MS) settings respectively. In (\RM) setting, the isoperimetric rations $\eta(A)$ and $\eta'(A)$ are already defined in the beginning, see~\eqref{bdry3}. We then consider {\[\rho(A):=\min_{\substack{B\in\A\\B\subseteq A}} \eta(B)\,, \qquad\rho'(A):=\min_{\substack{B'\in\A\\B'\subseteq A}} \eta'(B')\]}
in (\FS) and (\MS)  settings. And in (\RM) setting we take
\[\rho(A):=\inf_{\substack{B\in\A\\ B\subset A\\ \bar B\cap \partial_iA=\emptyset}} \eta(B)\,,\qquad\rho'(A):=\inf_{\substack{B'\in\A\\ B'\subset A\\ \bar B'\cap \partial_iA=\emptyset}} \eta'(B')\]
The constant $\rho(A)$ in (\RM) setting  is the Cheeger constant of $A$ when the Dirichlet boundary condition on $\partial_iA$ is imposed, we refer to \cite{bus80,MR0397619} for more information on  the Cheeger constant on manifolds with Dirichlet and Neumann boundary conditions.
 We are now ready to define the \textit{higher order Cheeger--Steklov constants}.  For any $k\in\N$ and for any of three settings (\FS), (\MS) and (\RM),  we define the $k$-th Cheeger--Steklov constant of $M$ by 
\bq\iota_k(M)&\df&\inf_{(A_1,\cdots,A_k)\in\A_k}\max_{l\in\llbracket k\rrbracket}\rho(A_l)\rho'(A_l)\eq
where $\llbracket k\rrbracket:=\{1,\ldots, k\}$ and $\A_k$ is the set of all $k$-tuples $(A_1,\cdots,A_k)$ such that $\A_l\in \A$ for all $l\in\llbracket k\rrbracket$.
We recall the definition of the higher order Cheeger constants for the eigenvalues of a Markov generator in settings (\FS) and (\MS) and for the eigenvalues of the Laplace--Beltrami operator in setting (\RM):
\bq h_k(M)&\df&\inf_{(A_1,\cdots,A_k)\in\A_k}\max_{l\in\llbracket k\rrbracket}\eta(A_l)\eq
The sequence of the  higher order Cheeger constants is  called the \textit{connectivity spectrum}. One can see how closely $h_k$ and $\iota_k$ are related. We now state our main theorems. 

\begin{theorem}\label{A}
In setting (\FS),  there exists  a universal positive constant $c_0$ such that \bq
\fo k\in\lin v\rin,\qquad \sigma_k(M)&\geq&  \frac{c_0}{k^6} \frac{\iota_k(M)}{\lVe L\rVe}\eq
where $\lVe L\rVe$ is the largest absolute value of the elements of the diagonal of $L$.
\end{theorem}
The following theorem is an extension of Theorem \ref{A} to  setting (\MS).
\begin{theorem}\label{B}
In setting (\MS), there exists a universal positive constant  $c_1$ such that \bq
\fo k\in\N,\qquad \sigma_k(M)&\geq&  \frac{c_1}{k^6} \iota_k(M)\eq
\end{theorem}
The higher order Cheeger-Steklov inequality in setting (\RM) which is an extension of Escobar and Jammes results to higher Steklov eigenvalues states 
\begin{theorem}\label{C}
In setting (\RM), there exists a universal positive constant $c_2$ such that 
 \bq\fo  k\in \N,\qquad\sigma_k(M)&\ge&\frac{c_2}{k^6}\iota_k(M)\eq
\end{theorem}
 We recall that  for $k=2$, the Cheeger inequality in setting (\FS) was studied in \cite{MR875835,MR782626,MR743744}, and in settings~(\MS) in \cite{MR930082}, see also the lecture notes by Saloff-Coste \cite{MR1490046} for a review. The higher-order Cheeger inequality in setting (\FS) was conjectured by the second author \cite{MR2438701}, see also \cite{MR2914771}. This conjecture was proved by Lee,
Oveis Gharan and 
Trevisan~\cite{MR2961569}. Later, the second author \cite{MR3323580} extended their result to (\MS) and (\RM) settings; see also \cite{Fun} for the result on closed manifolds. The higher order Cheeger inequality in (\FS) setting for the operator $L$ states (see \cite[Theorem 3.8]{MR2961569} and \cite[Theorem 2]{MR3323580})
\bqn{chfs}
\fo  k\in\lin v\rin,\qquad \lambda_k(M)&\geq&  \frac{c_3}{{k^8}} \frac{h_k^2(M)}{{\lVe L\rVe}}\eqn
and in  (\MS) and (\RM)  settings states \cite{MR3323580}
\bqn{chms}
\fo k\in\N,\qquad \lambda_k(M)&\geq&  \frac{c_4}{k^6} h_k^2(M)\eqn
where $c_3$ and $c_4$ are universal positive constants.
As we mentioned before, our main results, Theorems \ref{A}, \ref{B} and \ref{C} for Steklov eigenvalues,  can be viewed as a counterpart of the higher order Cheeger inequalities for the Laplace spectrum.
 We remark that even for $k=2$, Theorem \ref{A} and Theorem \ref{B}   are new. \\
We now discuss about an improvement  of the dependency on $k$ in Theorems \ref{A}, \ref{B}, and \ref{C}. In \cite[Theorem 4.1]{MR2961569} and \cite[Theorem 13]{MR3323580}, it is shown that one can obtain a better lower bound when $\lambda_k$ is replaced by $\lambda_{2k}$ in \eqref{chfs} and \eqref{chms}
\bqn{2k-inq}\lambda_{2k}(M)&\ge&\begin{cases}
\frac{\tilde c_3}{{\log(k+1)}} \frac{h_k^2(M)}{{\lVe L\rVe}}&\quad\mbox{in setting (\FS)}\\
\frac{\tilde c_4}{\log^2(k+1)} h_k^2(M)&\quad\mbox{in settings (\MS) and (\RM)}
\end{cases}
\eqn
For Steklov eigenvalues we obtain analogous results.
\begin{prop}\label{proA} There are universal positive constants $\tilde c_1$ and $\tilde c_2$ such that
\bqn{2kinq}\sigma_{2k}(M)&\ge&\begin{cases}
\frac{\tilde c_1}{\log^2(k+1)} \frac{\iota_k(M)}{\lVe L\rVe}&\qquad\fo k\in\lin v\rin,\qquad\mbox{in setting (\FS)}\\
\frac{\tilde c_2}{\log^2(k+1)}\iota_k(M)&\qquad \fo k\in \N,\qquad\mbox{in settings (\MS) and (\RM)}
\end{cases}\eqn
\end{prop}
\begin{rem}
The sharpness of the coefficient of $h_{k}$ in (\ref{2k-inq}) was investigated in  \cite{MR3323580} using the noisy hypercube graph, and in \cite{MR2961569}  using the Ornstein--Uhlenbeck
process. Understanding the asymptotic sharpness of the coefficient of $\iota_{k}$ in \eqref{2kinq} is an interesting problem which needs a further investigation and remains open.  
\end{rem}
We now briefly discuss the idea of the proof of the main Theorems. To prove the main theorems we first introduce the Dirichlet-Steklov connectivity spectrum of $S$ on $M$. Second we show that eigenvalues of  $S$ can be viewed as a limit of  eigenvalues of  a family  of operators. Then we prove that the Dirichlet connectivity spectrum (introduced in \cite{MR2438701} and in \cite{MR3323580}) of this family of operators converges to Dirichlet-Steklov connectivity spectrum of $S$. Moreover, we show that this convergence is uniform in some sense. Then we  use the known lower bounds \cite{MR2961569,MR3323580} for eigenvalues of this family  of operators  in terms of their Dirichlet connectivity spectra to show that   the Steklov eigenvalues have similar lower bounds in term of the Dirichlet-Steklov connectivity spectrum. The final step is to relate the Dirichlet--Steklov connectivity spectrum to the higher order Cheeger--Steklov constants. This  is done using the co-area formula in each  setting (\FS), (\MS) and (\RM). Although the main idea of the proof in these three settings are the same,  the details and technicalities that we need to deal with in each setting are different. This makes the investigation of each setting interesting in its own and not only as a straightforward consequence of another setting. We aim to explore a deeper underlying connection between these three settings in future studies.\\

 It is also  interesting to study the higher order Cheeger-Steklov inequality when $L$  is a diffusion operator and when we also have a density on $V$. Here the associated Dirichlet--to--Neumann map $S$ (known  as the voltage--to--current map)  appears in the study of the electrical impedance tomography \cite{MR1955896,MR3460047}.   The techniques and methods that we develop in this paper can be used to obtain the higher order Cheeger--Steklov inequality  in this setting in terms of a weighted version of the higher order Cheeger--Steklov constants. The classical Cheeger inequality for weighted manifolds is studied in \cite{MR783536}, see also \cite{MR1646764,MR3323580}. We will address this in more details in a forthcoming work. \\
 
The paper is organized as follows. Section \ref{FS} deals with (\FS) setting and the proof of Theorem \ref{A} and Proposition \ref{proA}. In Section \ref{MS} we extends results in (\FS) setting to (\MS) setting. We also show that under the Dirichlet gap assumption on $M\setminus V$ the proof of Theorem \ref{B} can be simplified.   In Section \ref{RM} we prove Theorem \ref{C}. We also provide examples which show the necessity of both isoperimetric ratios appearing in the definition of $\iota_k$.    Although the ideas and techniques in three sections  \ref{FS},  \ref{MS}, and  \ref{RM} are related, the reader does not need to read the sections in order. 
In Appendix \ref{App}, we discussed some continuity properties in setting (\MS). In particular, it shows that in this setting under somewhat  restrictive conditions, a Steklov operator  without the Dirichlet gap assumption on $M\setminus V$ can be viewed as a limit of  a family of Steklov operators with the Dirichlet gap assumption on $M\setminus V$. 

\subsection*{Acknowledgement} 
The authors are grateful to the organizers of  the program  ``Interactions between Partial Differential Equations \& Functional Inequalities" in the Mittag--Leffler Institute, and to the organizers of the conference ``Curvature-dimension in Lyon~1" in Institut Camille Jordan at Université Claude Bernard Lyon~1, as well as  the hospitality of  these two institutes  where main part of the research was conducted. The authors thank Luigi Provenzano for providing reference \cite{Pthese} for a detailed proof of Theorem \ref{provenzano}.
The first named author is supported by the Mittag--Leffler institute through the European Postdoctral Fellowship. She gratefully acknowledges  its support. She also thanks the Max Planck Institute for Mathematics in Bonn for its support during the starting phase of the project.  The second author is supported by the ANR STAB (Stabilité du comportement asymptotique d'EDP, de processus stochastiques et de leurs discrétisations : 12-BS01-0019).

\section{The finite state space framework}\label{FS}

Let $L\df L(x,y)_{x,y\in M}$ be an irreducible Markov generator on the finite set $M$. Recall that $L$ is Markovian if 
$$\fo x\neq y\in M\quad L(x,y)\ge0,\qquad\hbox{and} \, \qquad \sum_{y\in M} L(x,y)=0$$
and is called irreducible if for every $x,y\in M$ there exists a sequence $x=x_0,x_1,\ldots,x_l=y$ of elements of $M$ such that $L(x_j,x_{j+1})>0$ for any $j\in\llbracket0,l-1\rrbracket:=\{0,\ldots,l-1\}$.
 Denote by $\mu\df(\mu(x))_{x\in M}$ its unique invariant probability,
characterized by
\bq
\fo y\in M,\qquad \sum_{x\in M}\mu(x) L(x,y)&=&0\eq
Let $V$ be a proper subset of $M$, i.e.\ $\emptyset \varsubsetneq V\varsubsetneq M$.
Define the corresponding Steklov operator $S$ on $\cF(V)$, the space of functions on $V$, via the following procedure.
Given $f\in \cF(V)$, let $F$ be its harmonic extension on $M$, namely the unique $F\in \cF(M)$ satisfying
\bqn{F}
\begin{cases}
L[F](x)=0\,,&\,\hbox{if}\; \; x\in M\setminus V\\
F(x)=f(x)\,,&\,\hbox{if}\; \; x\in V
\end{cases}
\eqn
Then we consider
\bqn{defFS}
\fo x\in V,\qquad S[f](x)&\df& L[F](x)\eqn
The following observation should be classical.
\begin{pro}\label{pro1}
The operator $S$ is an irreducible Markov generator on $V$ whose invariant measure
is $\nu$, the renormalized restriction of $\mu$ to $V$. 
\end{pro}
\par
Assume that  $\mu$ is furthermore reversible for $L$, namely
\bq
\fo x,y\in M,\qquad \mu(x)L(x,y)&=&\mu(y)L(y,x)\eq
It follows that $S$ is equally reversible with respect to $\nu$, and the spectra of  $-S$  and $-L$ are non-negative.
Denote by $0=\sigma_1,\sigma_2, \sigma_3, \dots, \sigma_v$, with $v\df\card(V)$, the  eigenvalues of $-S$ in $\RR$ with multiplicities,
indexed 
so that $0=\sigma_1<\sigma_2\leq \sigma_3\leq \cdots \leq \sigma_v$. 
\par
Our goal is to investigate these eigenvalues.
Follows a way to approximate them.\par
For any $r>0$, consider the Markov generator defined by
\bq
\fo x\neq y\in M,\qquad L\rr(x,y)&\df& 
\begin{cases}
rL(x,y)\,,&\hbox{ if $x\in M\setminus V$}\\
L(x,y)\, ,&\hbox{ if $x\in V$}
\end{cases}\eq
Since $\mu$ is reversible for $L$, we will see (in Lemma \ref{lem1}) that $L^{(r)}$ is reversible with respect to its  invariant measure $\mu^{(r)}$. Hence the eigenvalues of $-L^{(r)}$ are non-negative. 
Let  $0=\lambda_1\rr,\lambda_2\rr, \lambda_3\rr, \dots, \lambda_m\rr$, with $m\df\card(M)$, be the  eigenvalues of $-L\rr$  in $\RR$ with multiplicities,
indexed 
so that $0=\lambda_1\rr<\lambda_2\rr\leq \lambda_3\rr\leq \cdots \leq \lambda_m\rr$.
\begin{pro}\label{pro2}
Assume that $L$ is reversible. For any $k\in \lin v\rin\df\{1, ..., v\}$, we have
\bq
\lim_{r\ri+\iy}\lambda_k\rr&=&\sigma_k\eq
and for any $k\in\lin m\rin\setminus\lin v\rin$,
\bq
\lim_{r\ri+\iy}\lambda_k^{(r)}&=&+\iy\eq
\end{pro}
\begin{rem}We believe that the above proposition should be true in the non-reversible case (where in the last convergence,
$\lambda_k^{(r)}$ is replaced by its real part). 
\end{rem}
We would like to estimate these eigenvalues via Cheeger type inequalities.
 Denote by $\cA$ the set of nonempty subsets from $M$.
We associate to any $A\in \cA$ a Dirichlet-Steklov operator $S_A$ on $\cF(A\cap V)$ in the following way:
given $f\in \cF(A\cap V)$, consider $F\in\cF(M)$ such that
\bqn{F2}
\begin{cases}
L[F](x)=0\,, &\quad\hbox{if}\;\, x\in A\setminus V\\
F(x)=0\,,&\quad\hbox{if}\;\, x\in M\setminus A\\
F(x)=f(x)\,,&\quad\hbox{if}\;\, x\in A\cap V\\
\end{cases}
\eqn
The existence and uniqueness of
 such a $F$ are similar to those of the solution of \eqref{F}, see e.g.\ the proof of Proposition \ref{pro1}.
 Indeed, one is brought back to this situation by replacing $V$ by  $V\cup (M\setminus A)$ and by extending $f$ to this set
 by making it vanish on $M\setminus A$.\par
Next define 
\bq
\fo x\in A\cap V,\qquad S_A[f](x)&\df& L[F](x)\eq
When $A\cap V\neq\emptyset$, we will check that $S_A$ is always a subMarkovian generator (i.e. $S_A(x,y)\ge0$, for any  $x\neq y$, and $\sum_{y\in V} S_A(x,y)\le0$) maybe not irreducible, but 
Perron-Frobenius' theorem  enables to consider the smallest eigenvalue $\sigma_1(A)$ of $-S_A$. By convention, when $A\cap V=\emptyset$, $\cF(\emptyset)\df \{0\}$ and $\sigma_1(A)=+\infty$.
Next we introduce the Dirichlet--Steklov connectivity spectrum $(\kappa_1, \kappa_2, ..., \kappa_v)$ of $S$ via
\bqn{kappak}
\fo k\in\lin v\rin,\qquad 
\kappa_k&\df&\min_{(A_1, ..., A_k)\in\cA_k}\max_{l\in\lin k\rin} \sigma_1(A_l)\eqn
where $\cA_k$ is the set of $k$-tuples $(A_1, A_2, ..., A_k)$ of disjoints elements from $\cA$. Notice that definition \eqref{kappak} can be written as 
\bqn{kappak-a}
\fo k\in\lin v\rin,\qquad 
\kappa_k&\df&\min_{(A_1, ..., A_k)\in\cA_k(V)}\max_{l\in\lin k\rin} \sigma_1(A_l)\eqn
where $\cA_k(V)$ is the set of all disjoint $k$-tuple in 
$
\cA(V)\df \{A\in \cA\St A\cap V\in\cA\}$.
The above definitions are valid in all generality, but (for the moment) they are mainly useful
under the reversibility assumption:
\begin{theo}\label{theo1}
Assume that $L$ is reversible.
There exists a universal constant $c>0$ such that 
\bq
\fo k\in\lin v\rin,\qquad  \frac{c}{k^6} \kappa_k\le\sigma_k \le\kappa_k\eq
\end{theo}
\par
The interest of the Dirichlet--Steklov connectivity spectrum is that it is strongly related to higher order inequalities.
We need further definitions.
Introduce
the boundary of any $A\in\cA$ via
\bq
\pa A&\df& \{(x,y)\St x\in A, \, y\in A\cc\}\eq
Consider the measure $\umu$ defined on $M\times M$ by
\bqn{bdry1}
\fo x,y\in M,\qquad
\umu(x,y)&=& \lt\{
\begin{array}{ll}
\mu(x)L(x,y)\,,&\hbox{ if $x\neq y$}\\
0\,,&\hbox{ if $x=y$}
\end{array}\rt.
\eqn
it enables to measure $\pa A$ through $\umu(\pa A)$. As a consequence, we can define the isoperimetric ratios
\bq
\eta(A)&\df& \frac{\umu(\pa A)}{\mu( A)}\\
\eta'(A)&\df& \frac{\umu(\pa A)}{\mu(  A\cap V)}
\eq
By convention  $\eta'(A)=+\iy$ if $A\cap V=\emptyset$.
The ratio $\eta'(A)$ is the discrete analogue of quantities introduced by Escobar \cite{MR1479552} and Jammes \cite{MR3449183},
since in their terminology, $\pa A$ and $A\cap V$ can be seen respectively as the interior and  exterior boundaries,
when the set $V$ itself is seen as a boundary of $M$.
\par
Next consider
\[\rho(A):=\min_{\substack{B\in\A\\B\subseteq A}} \eta(B)\]\[\rho'(A):=\min_{\substack{B'\in\A\\B'\subseteq A}} \eta'(B')\]

\par
For any $k\in\lin v\rin$, introduce the $k$-th Cheeger--Steklov constant of $V$ by
\bq
\iota_k&\df&\min_{(A_1, ..., A_k)\in\cA_k}\max_{l\in\lin k\rin} \rho(A_l)\rho'(A_l)\eq
Remark that $\iota_1=0$ by taking $A=M$. The next result can be seen as an extension to higher order Cheeger inequalities (in the discrete case) of Théorème 1 of Jammes \cite{MR3449183}:
\begin{theo}\label{theo2}
Assume that $L$ is reversible and let $c$ be the constant of Theorem \ref{theo1}.
We have 
\bq
\fo k\in\lin v\rin,\qquad \sigma_k&\geq & \frac{c}{k^6} \frac{\iota_k}{\lVe L\rVe}\eq
where $\lVe L\rVe$ is the largest absolute value of the elements of the diagonal of $L$.
\end{theo}
Let consider 
\bq
h_k'&\df&\min_{(A_1, ..., A_k)\in\cA_k(V)}\max_{l\in\lin k\rin} \eta'(A_l)\eq
\begin{pro}\label{hkup}
Assume that $L$ is reversible.
We have 
\bq
\fo k\in\lin v\rin,\qquad \sigma_k&\le & h'_k\eq
\end{pro}

\begin{rem}
Let $L$ be a reversible Markov generator but not necessarily irreducible. Let $X\df(X_t)_{t\geq 0}$ be a Markov process generated by $L$, starting from $x$ under the probability $\PP_x$.
Assume that the reaching time of $V$ denoted by $\tau$:
\bq
\tau&\df& \inf\{t\geq 0\St X_t\in V\}\eq  is almost surely finite. Then all of the results above are valid without irreducibility condition. In particular, $\sigma_k=0$ if and only if $h'_k=0$. Indeed one way is obvious  due to Proposition \ref{hkup}.  For the ``only if'' part,   $\sigma_k=0$, implies $\iota_k=0$ by Theorem \ref{theo2}. Therefore there exists $(A_1, ..., A_k)\in\cA_k(V)$ such that $\umu(\partial A_l)=0$ for all $l\in\lin k\rin$. It follows $h'_k=0$.    Note that the number of zeros determines the number of communicating classes. Recall that for the eigenvalues of $L=L^{(1)}$, the result of Lee,
Oveis Gharan and Trevisan~\cite{MR2961569} implies that $\lambda_k=0$ if and only if the $k$-th Cheeger constant $h_k$$$h_k:=\min_{(A_1, ..., A_k)\in\cA_k}\max_{l\in\lin k\rin} \eta(A_l)$$ is zero. In comparison, we see that the $h'_k$ plays the role of $h_k$ for the Steklov problem . \end{rem}
\prooff{Proof of Proposition \ref{pro1}}
It is based on the following simple probabilistic interpretation of $S$.
Let $X\df(X_t)_{t\geq 0}$ be a Markov process generated by $L$, starting from $x$ under the probability $\PP_x$.
Denote by $\tau$ its reaching time of $V$:
\bq
\tau&\df& \inf\{t\geq 0\St X_t\in V\}\eq
it is a.s.\ finite, since $L$ is irreducible.
A usual application of the martingale problem associated to $X$ shows that
for any function $G\in\cF(M)$, we have
\bq
\EE_x[G(X_\tau)]&=&G(x)+\EE_x\lt[\int_0^\tau L[G](X_s)\, ds\rt]\eq
In particular, for any  $f\in\cF(V)$, it appears that its harmonic extension defined in \eqref{F} is given by
\bq
\fo x\in M,\qquad F(x)&=&\EE_x[f(X_\tau)]\\
&=&\nu_x[f]\eq
where $\nu_x$ is the law of $X_\tau$ under $\PP_x$.
More precisely, we get the existence and uniqueness of the solution of \eqref{F},
even without assuming that $L$ is irreducible (only the finiteness of $\tau$ is needed).
We deduce that for any $f\in\cF(V)$ and any $x\in V$,
\bq
S[f](x)&=&\sum_{y\in M\setminus\{x\}} L(x,y)(F(y)-F(x))\\
&=&\sum_{y\in M\setminus\{x\}}\sum_{z\in V} L(x,y) \nu_y(z)(f(z)-f(x))\eq
namely, the matrix associated to $S$ is given by
\bq
\fo x,z\in V,\qquad S(x,z)&\df& \lt\{
\begin{array}{ll}
\sum_{y\in M\setminus\{x\}} L(x,y) \nu_y(z)\,,&\hbox{ if $x\neq  z$}\\
-\sum_{y\in V\setminus\{x\}}S(x,y)\,,&\hbox{ if $x=z$}\end{array}\rt.\eq
On this expression, it is clear that $S$ is a Markov generator, namely that it satisfies $S(x,z)\geq 0$ for any
$x\neq z\in V$ and  $\sum_{z\in V}S(x,z)=0$ for any $x\in V$.
It is also irreducible: for any $x,z\in V$, let $x_0=x, x_1, x_2, ..., x_l=z$ be a sequence of elements of $M$ such that
$L(x_j,x_{j+1})>0$ for any $j\in \lin 0, l-1\rin$. Let $(y_j)_{j\in \lin 0, k\rin}$ be the subsequence of $(x_j)_{j\in \lin 0, l\rin}$
consisting of the elements belonging to $V$. We have $y_0=x$, $y_k=z$ and from the above description of $S$,
it follows that $S(x_j,x_{j+1})>0$ for any $j\in \lin 0, k-1\rin$.
\par
It remains to check that $\nu$, the renormalized restriction of $\mu$ to $V$, is invariant for $S$.
For any $f\in \cF(V)$, we have, with $F$ constructed as in \eqref{F},
\bq
\nu[S[f]]&=&\frac1{\mu(V)}\sum_{x\in V}\mu(x)S[f](x)\\
&=&\frac1{\mu(V)}\sum_{x\in V}\mu(x)L[F](x)\\
&=&\frac1{\mu(V)}\sum_{x\in M}\mu(x)L[F](x)\\
&=&\frac{\mu[L[F]]}{\mu(V)}\\
&=&0\eq
It shows that $\nu$ is invariant for $S$.
\wwtbp
\par
\begin{remm}{(probabilist point of view)}\label{rem2}
A Markov process $Y\df(Y_t)_{t\geq 0}$ associated to the generator $S$ and starting from $x\in V$ can be obtained from a
Markov process $X\df(X_t)_{t\geq 0}$ associated to the generator $L$ and also starting from $x$, by erasing  its passages in $M\setminus V$.
More precisely, let $(\tau_n)_{n\in\ZZ_+}$ be the sequence of jump intertimes of $X$: 
\bq
\tau_0&\df& 0\\
\fo n\in\ZZ_+,\qquad
\tau_{l+1}&\df& \inf\{t\geq 0\St X_{ t+\tau_l}\neq X_{\tau_l}\}\eq
Let $(N_n)_{n\in\ZZ_+}$ be the sequence of integers for which $X_{\tau_1+\tau_2+\cdots+\tau_{N_n}}\in V$
and
consider
\bq
\fo n\in\ZZ_+,\qquad \uptau_n&\df& \sum_{p\in\lin n\rin} \tau_{N_p}\eq
Then we can construct the Markov process $Y$ through the relation
\bq
\fo t\geq 0,\qquad Y_t&\df& X_{\tau_1+\tau_2+\cdots+\tau_{N_n}}\,,
\qquad \hbox{ if\,\, $t\in[\uptau_n,\uptau_{n+1}[$}\eq
This observation inspired the introduction of the generators $L\rr$, for $r>0$: heuristically the generator of $Y$ is $L^{(\infty)}$,
namely $X$ is accelerated with an infinite speed in $M\setminus V$ and only its passages on $V$ remain. \par
The above probabilistic interpretation also enables to see directly that $S$ is irreducible and that the invariant measure $\nu$ of $S$ is just $\mu$ conditioned on $V$.
Indeed, for the latter assertion, by the ergodic theorem, we must have a.s.\ 
\bq
\fo y\in V,\qquad \nu(y)&=&\lim_{t\ri+\iy} \frac1t\int_0^t \un_{\{y\}}(Y_s)\, ds\eq
so it follows that for any $y,z\in V$,
\bq 
\frac{\nu(y)}{\nu(z)}&=&\lim_{t\ri+\iy}\frac{\int_0^t \un_{\{y\}}(Y_s)\, ds}{\int_0^t \un_{\{z\}}(Y_s)\, ds}\\
&=&\lim_{t\ri+\iy}\frac{\int_0^t \un_{\{y\}}(X_s)\, ds}{\int_0^t \un_{\{z\}}(X_s)\, ds}\\
&=&\frac{\mu(y)}{\mu(z)}\eq
\end{remm}
\par
\begin{remm}{(analytic point of view)}\label{rem10}
Recall that the Dirichlet form associated to $L$ (and $\mu$) is the bilinear form $\cE_L$ given by
\bq
\fo F,G\in \cF(M),\qquad \cE_L(F,G)&\df& -\int F L[G]\, d\mu\eq
It is symmetrical, if and only if $\mu$ is reversible with respect to $L$.\par
The 
carré du champ associated to $L$ is the bilinear functional $\Gamma_L$ defined by
\bqn{Gamma}
\fo F,G\in \cF(M),\,\fo x\in M,\qquad \Gamma_L[F,G](x)&\df& L[FG](x)-F(x)L[G](x)-G(x)L[F](x)\eqn
It is not difficult to compute more explicitly that
\bq
\fo F,G\in \cF(M),\,\fo x\in M,\qquad \Gamma_L[F,G](x)&\df& \sum_{y\in M}L(x,y) (F(y)-F(x))(G(y)-G(x))\eq
In particular, when $F=G$,  the r.h.s.\ looks like a weighted discrete gradient square, explaining the name carré du champ.
\par
From \eqref{Gamma}, we get that
\bq
\fo F,G\in \cF(M),\qquad \int \Gamma_L[F,G]\, d\mu&=&\cE_L(F,G)+\cE_L(G,F)\eq
and in particular 
\bq
\fo F\in \cF(M),\qquad \int \Gamma_L[F]\, d\mu&=&2\cE_L(F,F)\eq
where $\Gamma_L[F]$ stands for $\Gamma_L[F,F]$.
Furthermore, when $\mu$ is reversible with respect to $L$, we get
\bq
\fo F,G\in \cF(M),\qquad \int \Gamma_L[F,G]\, d\mu&=&2\cE_L(F,G)\eq
\par
These definitions are valid for any finite Markov generator $L$ and we can consider similarly $\cE_S$ and $\Gamma_S$.
For any $f,g\in\cF(V)$, let $F$ and $G$  be their harmonic extensions. 
It is clear that
\bqn{cELcES}
\cE_S(f,g)&=&\frac{\cE_L(F,G)}{\mu(V)}\eqn
and as a consequence, we have 
\bq
\int \Gamma_S[f,g]\, d\nu&=&\frac1{\mu(V)}\int \Gamma_L[F,G]\, d\mu\eq
which is an important relation in the analytical approach to the usual Steklov (or Dirichlet to Neumann) operators.\par
It follows immediately from \eqref{cELcES} that $\nu$ is reversible for $S$ when $\mu$ is assumed to be reversible for $L$.
\end{remm}
\par\sm
Since for any $r>0$, the generator $L\rr$ is irreducible, it admits a unique invariant probability $\mu\rr$.
\begin{lem}\label{lem1}
The probability measure $\mu\rr$ is given by
\bq
\fo x\in M,\qquad \mu\rr(x)&=&\lt\{
\begin{array}
{ll}
\frac{\mu(x)}{Z_r}\,,&\hbox{ if $x\in V$}\\
\frac{\mu(x)}{rZ_r}\,,&\hbox{ if $x\,  \in M\setminus V$}
\end{array}\rt.\eq
where $Z_r\df \mu(V)+(1-\mu(V))/r$ is the normalisation constant.\par
Furthermore, if $\mu$ is reversible for $L$, then $\mu\rr$ is reversible for $L\rr$.
\end{lem}
\proof
These are consequences of more general facts: assume that $H\in\cF(M)$ is positive: $H>0$.
Consider the operator $HL$ acting on $\cF(M)$ via
\bq
\fo F\in\cF(M),\, \fo x\in M,\qquad HL[F](x)&\df& H(x)L[F](x)\eq
It is an irreducible Markov generator. Let $(1/H)\cdot\mu$ be the positive measure admitting $1/H$ for density with respect to $\mu$.
We have
\bq
\fo F\in \cF(M),\qquad ((1/H)\cdot\mu)[HL[F]]&=&\mu[L[F]]\\
&=&0\eq
Thus the invariant probability measure of $HL$ is proportional to $(1/H)\cdot\mu$.
\par
Considering $H\df \un_V+r\un_{M\setminus V}$ (where $\un_V$ is the indicator function of $V$) leads to the first announced result.\par
For the second result, note that in general, when $\mu$ is reversible for $L$,
\bq
\fo F,G\in\cF(M),\qquad ((1/H)\cdot\mu)[F(HL)[G]]&=&\mu[FL[G]]\\
&=& \mu[GL[F]]\\
&=&((1/H)\cdot\mu)[G(HL)[F]]\eq
\wwtbp
\par\me
\prooff{Proof of Proposition \ref{pro2} }
In the reversible case, $-L$ is diagonalisable with real eigenvalues.
In view of Lemma \ref{lem1}, for any $r>0$, the same is true for $-L\rr$, denote by 
$0=\lambda_1\rr<\lambda_2\rr\leq \lambda_3\rr\leq \cdots \leq \lambda_m\rr$ its eigenvalues.
Let $\un=\Phi_1\rr, \Phi_2\rr, \Phi_3\rr, \ldots, \Phi\rr_m$ be corresponding eigenvectors. They are not unique (especially in the case of multiplicities larger than 1),
but we can and do choose them so that they are orthogonal with respect to $\mu\rr$:
\bq
\fo r\in(0,+\iy),\,\fo k\neq l\in\lin m\rin,\qquad 
\mu\rr[\Phi\rr_l\Phi\rr_k]&=&0\eq
Renormalize them with respect to the supremum norm $\lVe\cdot\rVe_\iy$ instead of the $\LL^2(\mu\rr)$ norm:
\bq
\fo r\in(0,+\iy),\,\fo l\in\lin m\rin,\qquad \lVe \Phi_l\rr\rVe_{\iy}&=&1\eq
\par
Consider $l\in \lin m\rin$ such that
\bqn{l}
\lt\{\begin{array}{rcl}
\liminf_{r\ri+\iy} \lambda_l\rr&<&+\iy\\
\liminf_{r\ri+\iy} \lambda_{l+1}\rr&=&+\iy
\end{array}\rt.
\eqn
By compactness, we can find an increasing sequence of positive numbers $(r_n)_{n\in\NN}$
and for any $k\in\lin l\rin$,
 a non-negative number $\lambda^{(\iy)}_k\in[0,+\iy)$ and a positive  function $\Phi^{(\iy)}_k\in\cF(M)$ with $\lVe \Phi^{(\iy)}_k\rVe_\iy=1$ such that
\bq
\lim_{n\ri\iy} r_n&=&+\iy\\
\lim_{n\ri\iy} \lambda^{(r_n)}_k&=&\lambda^{(\iy)}_k\\
\lim_{n\ri\iy} \Phi^{(r_n)}_k&=&\Phi^{(\iy)}_k\eq
Passing to the limit in the relations
\bq
\fo x\in V,\qquad L[\Phi^{(r_n)}_k](x)&=&L^{(r_n)}[\Phi^{(r_n)}_k](x)\\
&=&-\lambda^{(r_n)}_k\Phi^{(r_n)}_k(x)\eq
we get
\bq
\fo x\in V,\qquad L[\Phi^{(\iy)}_k](x)
&=&-\lambda^{(\iy)}_k\Phi^{(\iy)}_k(x)\eq
For $x\in M\setminus V$, we have instead
\bq
r_nL[\Phi^{(r_n)}_k](x)
&=&-\lambda^{(r_n)}_k\Phi^{(r_n)}_k(x)\eq
Since the r.h.s.\ converges to $-\lambda^{(\iy)}_k\Phi^{(\iy)}_k(x)$
for large $n\in\NN$, we deduce that
\bq
\fo x\in M\setminus V,\qquad L[\Phi^{(\iy)}_k](x)&=&\lim_{n\ri\iy} L[\Phi^{(r_n)}_k](x)\\
&=&0\eq
Thus denoting $\varphi_k$ the restriction of $\Phi^{(\iy)}_k$ to $V$, it appears that $\Phi^{(\iy)}_k$ is the harmonic extension of $\varphi_k$.
Note that $\varphi_k\neq 0$, otherwise we would conclude that $\Phi^{(\iy)}_k=0$, in contradiction with $\lVe \Phi^{(\iy)}_k\rVe_\iy=1$.
Thus $\lambda^{(\iy)}_k$ is an eigenvalue of $-S$.
Furthermore, passing to the limit in the relations
\bq
\fo j\neq  k\in\lin l\rin,\qquad \mu^{(r_n)}[\Phi^{(r_n)}_j\Phi^{(r_n)}_k]&=&0\eq
we see that
\bq
\fo j\neq  k\in\lin l\rin,\qquad \nu[\varphi_j\varphi_k]&=&0\eq
It follows that the $\lambda^{(\iy)}_k$, for $k\in\lin l\rin$, correspond to different eigenvalues of $-S$ (with multiplicities).
Namely, there exists an increasing mapping $N\St \lin l\rin\ri\lin v\rin$ (recall that $v\df \card(V)$) such that
\bq
\fo k\in\lin l\rin,\qquad \lambda_k^{(\iy)}&=&\sigma_{N(k)}\eq
and in particular, $v\geq l$. Conversely, consider
$\psi_1, \psi_2, ..., \psi_v$ a basis of $\cF(V)$ consisting of eigenvectors of $-S$ associated respectively to the eigenvalues
$\sigma_1, \sigma_2, ..., \sigma_v$. Since $\nu$ is reversible for $S$, we can and do choose these functions to be orthogonal in $\LL^2(\nu)$.
Let $\Psi_1, \Psi_2, ..., \Psi_v$ be the harmonic extensions of $\psi_1, \psi_2, ..., \psi_v$.
We furthermore impose that $\lVe \Psi_k\rVe_\iy=1$ for all $k\in\lin v\rin$.
Consider the vector space $W\subset \cF(M)$ generated by these functions
\bq
W&\df& \Vect( \Psi_k\St k\in\lin v\rin)\eq
 Due to the variational principle, we have for any $r>0$,
 \bq
 \lambda_v\rr&\leq & \sup_{F\in W\setminus\{0\}}\frac{-\mu\rr[FL\rr[F]]}{\mu\rr[F^2]}\eq
 Since the functions from $W$ are harmonic on $M\setminus V$, we have for any $r>0$,  with the notation of Lemma~\ref{lem1},
 \bq
 \fo F\in W,\qquad
 -\mu\rr[FL[F]]&=&-\frac{\mu(V)}{Z_r}\nu[FL[F]]\\
 &=&-\frac{\mu(V)}{Z_r}\nu[fS[f]]\\
 &\leq & \frac{\mu(V)}{Z_r}\sigma_v\nu[f^2]\eq
 where $f$ is the restriction of $F$ to $V$. We also have
 \bq
 \mu\rr[F^2]&=&\frac{ \mu[\un_V f^2]+\mu[\un_{M\setminus V} F^2]/r}{Z_r}\\
 &\geq & \frac{\mu(V)}{Z_r}\nu[f^2]\eq
 We deduce from these two bounds that
 \bq
 \lambda\rr_v&\leq & \sigma_v\eq
 and 
 \bqn{lsi}
\limsup_{r\ri+\iy} \lambda_v\rr&<&+\iy\eqn
i.e.\ $l\geq v$
and finally $l=v$.
\par
It follows that 
\bqn{conv}
\fo k\in\lin v\rin,\qquad \lim_{n\ri\iy} \lambda^{(r_n)}_k&=&\sigma_k\eqn
Taking into account \eqref{lsi}, for any increasing subsequence $(R_n)_{n\in\NN}$ of positive numbers diverging to $+\iy$, we can extract 
another subsequence $(r_n)_{n\in\NN}$ such that \eqref{conv} is true, we conclude by compactness that
\bq
\fo k\in\lin v\rin,\qquad \lim_{r\ri+\iy} \lambda^{(r)}_k&=&\sigma_k\eq
The last assertion of Proposition \ref{pro2} is a consequence of $l=v$ and of the definition of $l$ in \eqref{l}.
\wwtbp
\par\me
Before coming to the proof of Theorem \ref{theo1}, let us check that for any $A\in\cA(V)$, $S_A$ is a subMarkovian generator.
The argument is similar to that of the proof of Proposition \ref{pro1} and is based on the probabilistic representation of the solution $F$
of \eqref{F2}:
\bqn{solF2}
\fo x\in M,\qquad F(x)&=&\EE_x[f(X_{\tau_{A\cap V}})\un_{\tau_{A\cap V}<\tau_{M\setminus A}}]\eqn
where $(X_t)_{t\geq 0}$ is a Markov process generated by $L$ and starting from $x$, and
for any $B\subset M$, $\tau_B$ is the hitting time of $B$:
\bq
\tau_B&\df& \inf\{t\geq 0\St X_t\in B\}\eq
As a consequence, the first eigenvalue $\sigma_1(A)$ of $-S_A$ is non-negative.
It vanishes, if and only if there is no path (whose transitions are permitted by $L$) going out of $A$ without passing through $A\cap V$.
\par
Assume that $\mu$ is reversible with respect to $L$. By the variational formulation of eigenvalues and using the notation of Remark \ref{rem2},
we have for $A\in\cA$,
\bqn{sigma1}
\sigma_1(A)&=&\inf\lt\{\frac{\cE_{S_A}(f,f)}{\nu_{A\cap V}[f^2]}\St f\in \cF(A\cap V)\rt\}\eqn
where $\nu_{A\cap V}$ is the renormalized restriction of $\mu$ to $A\cap V$, which is reversible with respect to $S_A$.
As in \eqref{cELcES}, in the above formula, $\cE_{S_A}(f,f)$ can be replaced by $\cE_L(F,F)/\mu(A\cap V)$, where $F$ is associated to 
$f$ via \eqref{F2}.
\par
We can now come to the
\prooff{Proof of Theorem \ref{theo1}}
The upper bound of $\sigma_k$ is  a direct consequence of the variational characterization of $\sigma_k$
\bq
\sigma_k&=&\min_{H\in\cF_k(V)}\max_{f\in H\setminus\{0\}}\frac{\cE_{S}(f,f)}{\nu[f^2]}\eq
where $\cF_k(V)$ is the set of all $k$-dimensional subspace of $\cF(V)$, by taking $H$ as the space spanned by the first eigenfunctions of $S_{A_l}$, $l\in\llbracket k\rrbracket$.\\
The proof of the lower bound is based on the higher order Dirichlet-Cheeger inequalities for finite irreducible and reversible Markov generators.
So assume that $\mu$ is reversible with respect to $L$
and let  $0=\lambda_1(L)<\lambda_2(L)\leq \lambda_3(L)\leq \cdots \leq \lambda_m(L)$ be the eigenvalues of $-L$.
Associate to any $A\in\cA$ its first Dirichlet eigenvalue
\bq
\lambda_1(A)&\df& \inf\lt\{\frac{\cE_{L}(F,F)}{\mu[F^2]}\St F\in \cF(M)\hbox{ with }F\hbox{ vanishing on }M\setminus A\rt\}\eq
This is the same definition as \eqref{sigma1} if we had taken $V=M$.
Next define for any $k\in\lin m\rin$,
\bq
\Lambda_k(L)&\df&\min_{(A_1, ..., A_k)\in\cA_k}\max_{l\in\lin k\rin} \lambda_1(A_l)\eq
The higher order Dirichlet-Cheeger inequalities of Lee, Gharan and Trevisan \cite{MR2961569} (see also \cite{MR3323580} for its Markovian reformulation) assert that there exists a universal constant $c>0$ such that
\bq
\fo k\in\lin m\rin,\qquad \lambda_k(L)&\geq & \frac{c}{k^6} {\Lambda_k(L)}\eq
In particular, we can apply them to $L\rr$ for $r>0$:
\bqn{DCr}
\fo k\in\lin m\rin,\qquad \lambda^{(r)}_k\ =\ \lambda_k(L\rr)\ \geq \ \frac{c}{k^6} \Lambda_k(L\rr)=:\Lambda^{(r)}_k\eqn
From Proposition \ref{pro2}, we know the behavior for large $r>0$ of the l.h.s., for $k\in\lin v\rin$, so
it remains to investigate the r.h.s.\par
Fix $A\in\cA$ and consider for $r>0$,
\bq
\lambda_1\rr(A)&\df& \inf\lt\{\frac{\cE_{L\rr}(F,F)}{\mu\rr[F^2]}\St F\in \cF(M)\hbox{ with }F\hbox{ vanishing on }M\setminus A\rt\}\eq
It is the smallest eigenvalue of $-L\rr_A$, where $L\rr_A$ is the subMarkovian generator acting on $\cF(A)$ whose matrix is the $(A\times A)$-restriction  of the matrix corresponding to $L\rr$.
The proof of Proposition \ref{pro2} can easily be adapted to this situation to show that as $r$ goes to $+\iy$, the first $\card(A\cap V)$  eigenvalues of $-L\rr_A$
converge to the eigenvalues of $-S_A$. In particular we get
\bq
\lim_{r\ri+\iy} \lambda_1\rr(A)&=&\sigma_1(A)\eq
Since $\cA_k$ is a finite set, it follows that
\bq
\fo k\in \lin v\rin,\qquad \lim_{r\ri+\iy} \Lambda^{(r)}_k&=& \kappa_k\eq
where the r.h.s.\ is defined in \eqref{kappak}. The wanted result is thus obtained by passing to the limit in \eqref{DCr} as $r$ goes to $+\iy$.
\wwtbp
\par
\prooff{Proof of Theorem \ref{theo2}}
To relate the $\kappa_k$, for $k\in\lin v\rin$, to  isoperimetric quantities, we will adapt a computation of Jammes \cite{MR3449183} to the finite setting.
Fix $A\in\cA$ and let us come back to \eqref{sigma1}.
More precisely, consider $f\in\cF(A\cap V)$ a minimizer of the infimum in the r.h.s.\ of \eqref{sigma1} and $F$  the associated solution of \eqref{F2}.
From the Perron-Frobenius' theorem, we know that we can and do  choose $f$ to be non-negative and from \eqref{solF2}, we also have $F\geq 0$.
We are looking for a lower bound on the ratio
\bq\frac{\cE_L(F,F)}{\mu[f^2\un_{A\cap V}]}&=&\frac{\sum_{x\neq y\in M}\mu(x)L(x,y)(F(y)-F(x))^2}{2\sum_{x\in A\cap V} \mu(x)f^2(x)}\eq
So multiply the numerator and the denominator by $\sum_{x'\neq y'\in M}\mu(x')L(x',y')(F(y')+F(x'))^2$.
In the numerator we get
\bqn{linkfm2}
\nonumber\lefteqn{\hskip-20mm \sum_{x'\neq y'\in M}\mu(x')L(x',y')(F(y')+F(x'))^2\sum_{x\neq y\in M}\mu(x)L(x,y)(F(y)-F(x))^2}\\
&\geq & \lt(\sum_{x\neq y\in M}\mu(x)L(x,y)(F(y)+F(x))\vert F(y)-F(x)\vert\rt)^2\\
\nonumber&=&\lt(\sum_{x\neq y\in M}\mu(x)L(x,y)\vert F^2(y)-F^2(x)\vert \rt)^2\eqn
where for the first bound we used the Cauchy-Schwarz inequality with respect to the measure $\umu$ outside the diagonal of $M\times M$.
Concerning the denominator, we begin by noting that
\bqn{linkfm}
\nonumber\sum_{x'\neq y'\in M}\mu(x')L(x',y')(F(y')+F(x'))^2&\leq &2\sum_{x'\neq y'\in M}\mu(x')L(x',y')(F^2(y')+F^2(x'))\\
\nonumber&=&4\sum_{x'\neq y'\in M}\mu(x')L(x',y')F^2(x')\\
\nonumber&=&4\sum_{x'\in M}\mu(x') \lve L(x',x')\rve F^2(x')\\
&\leq &4\lVe L\rVe \sum_{x'\in M}\mu(x')F^2(x')\eqn
where we used the reversibility of $\mu$ with respect to $L$ for the first equality.
For any $G\in\cF(M)$, denote $\vert d G\vert $ the function on $M\times M$ given by
\bq
\fo (x,y)\in M,\qquad\vert dG\vert(x,y)&\df& \vert G(y)-G(x)\vert \eq
Putting together the above computations, we have obtained
\bq
\sigma_1(A)&\geq &\frac1{8\lVe L\rVe}\frac{\umu[\vert dF^2\vert]}{\mu[F^2]}\frac{\umu[\vert dF^2\vert]}{\mu[f^2\un_{A\cap V}]}\eq
To deal with the ratios of the r.h.s., recall the co-area formula (see for instance Formula (3.3.2) page  381 of the lecture notes of Saloff-Coste \cite{MR1490046}):
for any non-negative $G\in \cF(M)$ vanishing somewhere, we have
\bq
\umu[\vert dG\vert]&=&\int_0^{\tau} \umu [\pa D_t]\, dt\eq
where
\bqn{linkfm3}
\nonumber\fo t\geq 0,\qquad D_t&\df& \{x\in M\St G(x)\geq t\}\\
\tau&\df& \inf\{t\geq 0\St D_t =\emptyset\}\\
\label{linkfm4}&=&\inf\{t>0\St \umu(\pa D_t)=0\}
\eqn
We also have \bq
\mu[G]&=&\int_0^{\tau} \mu[D_t]\, dt\eq
Applying these formulas with $G\df F^2$ (which vanishes somewhere since $A\neq M$), we deduce that
\bq
\frac{\umu[\vert dF^2\vert}{\mu[F^2]}&\geq & \inf\lt\{ \frac{ \umu(\pa D_t)}{\mu[D_t]}\St t\geq 0\rt\}\\
&\geq & \min\lt\{\eta(B)\St B\in\cA,\,B\subset A\rt\}\eq
since we have $D_t\subset A$ for all $t\geq 0$.
Furthermore we  have
\bq
\mu[f^2\un_{A\cap V}]&=&\mu[F^2\un_{A\cap V}]\\
&=&\int_0^{+\iy} \mu[D_t\cap A\cap V]\, dt\\
&=&\int_0^{+\iy} \mu[D_t\cap V]\, dt\eq
so we deduce similarly that
\bq
\frac{\umu[\vert dF^2\vert]}{\mu[f^2\un_{A\cap V}]}&\geq & \min\lt\{\eta'(B)\St B\in\cA,\,B\subset A\rt\}\eq
Finally we have shown that
\bq
\fo A\in\cA,\qquad \sigma_1(A)&\geq & \frac{\rho(A)\rho'(A)}{8\lVe L\rVe}\eq
It follows that 
\bqn{aah}
\fo k\in\lin v\rin,\qquad \kappa_k&\geq & \frac{\iota_k}{8\lVe L\rVe}\eqn
and Theorem \ref{theo2} is now an immediate consequence of Theorem \ref{theo1}.\wwtbp
\prooff{Proof of Proposition \ref{hkup}}
Consider the variation characterization of $\sigma_k$:
\bq
\sigma_k&=&\min_{H\in\cF_k(V)}\max_{f\in H\setminus\{0\}}\frac{\cE_{S}(f,f)}{\nu[f^2]}=\min_{H\in\cF_k(V)}\max_{f\in H\setminus\{0\}}\frac{\cE_L(F_f,F_f)}{\mu[f^2\un_{ V}]}\eq
where $\cF_k(\cdot)$ is the set of all $k$-dimensional subspace of $\cF(\cdot)$, and $F_f$ is solution to \eqref{F2}, the harmonic extension of $f$ to $M\setminus V$. We can rewrite the variational characterisation in the following equivalent way. 
\bq
\sigma_k&=&\min_{\substack{H\in\cF_k(M)\\H|_{V}\in\cF_k(V)}}\max_{F\in H\setminus\{0\}}\frac{\cE_{L}(F,F)}{\mu[F^2\un_{V}]}\eq

Indeed for every $f\in \cF(V)$, and all $F\in \cF(M)$ with $F|_V=f$ we have
$$\cE_L(F_f,F_f)\le \cE_{L}(F,F)$$
 This is due to harmonic property of $F_f$, for a more detail see  \eqref{harmonicmin}. Let $(A_1, ..., A_k)\in\cA_k(V)$ and consider $H:=\Vect(\un_{A_l}\St l\in\llbracket k\rrbracket)\in\cF_k(M)$. It is also clear that $H|_V\in\cF_k(V)$.
\bq\frac{\cE_L(\un_{A_l},\un_{A_l})}{\mu[\un_{A_l\cap V}]}&=&\frac{\sum_{x\neq y\in M}\mu(x)L(x,y)(\un_{A_l}(y)-\un_{A_l}(x))^2}{2\mu(A_l\cap V)}\\&=&\frac{\sum_{x\in A_l,\, y\in A_l^c}\mu(x)L(x,y)+\mu(y)L(y,x)}{2\mu(A_l\cap V)}\\&=&\eta'(A_l)\eq
It implies 
\bq\sigma_k&\le& \min_{(A_1, ..., A_k)\in\cA_k(V)}\max_{l\in\lin k\rin} \eta'(A_l) =h_k'\eq and completes the proof. 
\wwtbp

We conclude this section by the proof  of        Proposition \ref{proA} in the introduction. 
\begin{pro}\label{fs-1} There is a universal positive constant $c'$  such that
\bq\fo k\in\lin v\rin,\qquad \sigma_{2k}&\ge&
\frac{C_1}{\log^2(k+1)} \frac{\iota_k}{\lVe L\rVe}\eq
\end{pro}
\proof By \cite[Theorem 4.6]{MR2961569} and \cite[Section 2]{MR3323580},  we have
\bq\fo k\in\lin v\rin,\qquad\lambda_k^{(r)}&\ge&\frac{c}{\log^2(k+1)}\Lambda^{(r)}_k\eq
where $c$ is a universal positive constant. Passing to limit and using \eqref{aah} we get 
\[\fo k\in\lin v\rin,\qquad\sigma_k=\lim_{r\to\infty}\lambda_k^{(r)}\ge\frac{c}{\log^2(k+1)}\kappa_k\ge\frac{c}{8\log^2(k+1)} \frac{\iota_k}{\lVe L\rVe}\]
and the statement follows.

\wwtbp\par
\section{The measurable state space framework}\label{MS}

Let $(M,\cM,\mu)$ be a probability measure space, endowed with a Markov kernel  $P$ 
 leaving $\mu$ invariant~(i.e. $\mu[P[F]]=\mu[F]$, for any bounded measurable function $F$). The Markov kernel $P$ defines a map $P:\LL^2(\mu)\to \LL^2(\mu)$ by $P[F](x):=\int_{M}P(x,dy)F(y)$. It has has the following  properties 
 \[P[\un]=\un,\quad\hbox{and}\quad \fo F\ge0\,\Rightarrow\, P[F]\ge0\]
 
We assume that $P$ is \textbf{weakly mixing}, in the following sense.
Let $Z\df (Z(n))_{n\in\ZZ_+}$ be a Markov chain whose transition kernel is $P$. As usual, we indicate that $Z$ is starting from $x\in M$, i.e.\ $Z(0)=x$, by putting $x$ in index of the underlying probability $\PP_x$ and expectation $\EE_x$ (more generally, this index will stand for the initial law of $Z(0)$.
Denote by $\cA$ the set of $A\in \cM$  such that $0<\mu(A)\le1$.
For any  $A\in\cA$, define
the \textbf{hitting time of $A$} by $Z$ via
\bqn{tauA}
\tau_{A}&\df&\inf\{n\in\ZZ_+\St Z(n)\in A\}\eqn
The weak mixing assumption asks for $\tau_{A}$ to be $\PP_x$-a.s.\ finite, for any $x\in M$ and any $A\in\cA$ (but what follows can be adapted to the situation where $\tau_{A}$ is a.s.\ finite, $\mu$-a.s.\ in $x\in M$ and for any $A\in\cA$).\par
Fix some $V\in \cA$, we introduce  corresponding \textbf{Steklov Markov kernel} $K$ and \textbf{Steklov generator} $S$ in the following way:
let $\cB(V)$ be the set of bounded measurable mappings defined on $V$.
To  any $f\in\cB(V)$, we associate the mapping $F_f\in\cB(M)$ given by
\bqn{Fgen}
\fo x\in M,\qquad F_f(x)&\df& \EE_x[f(Z(\tau_{V}))]\eqn
and we define
\bqn{Skern}
\fo x\in V,\qquad \lt\{
\begin{array}{rcl}K[f](x)&\df& P[F_f](x)\\
S[f](x)&\df& K[f](x)-f(x)
\end{array}\rt.\eqn

\par
Note that $K$ is a \textbf{Markov transition operator}, in the sense that it preserves the non-negativity of functions, as well as $\un_V$ (the mapping always taking the value 1 on $V$).
It is immediate to check that 
the function $F_f$ defined in \eqref{Fgen} is given by
\bq
F_f&=&\sum_{n\in\ZZ_+} (\un_{M\setminus V}P)^n\un_V[f]\eq
where the indicator functions are seen as multiplication operators. It follows that the
 transition kernel of $K$ is $\sum_{n\in\ZZ_+} (P\un_{M\setminus V})^nP\un_V$.
 The function $F_f$ is called the \textbf{harmonic extension} of $f$ to $M$, because we have
 \bqn{harm}
 \fo x\in M\setminus V,\qquad (P-I)[F_f](x)&=&0\eqn
 where $I$ stands for the identity operator (it will always be so in the sequel, even when the underlying space will not be the same). Indeed, we have
 on $M\setminus V$,
 \bq
 P[F_f]&=&\un_{M\setminus V}P[F_f]\\
 &=&\un_{M\setminus V}P\sum_{n\in\ZZ_+} (\un_{M\setminus V}P)^n\un_V[f]\\
 &=&\sum_{n\in\NN} (\un_{M\setminus V}P)^n\un_V[f]\\
 &=&\sum_{n\in\ZZ_+} (\un_{M\setminus V}P)^n\un_V[f]-\un_V[f]\\
 &=&\sum_{n\in\ZZ_+} (\un_{M\setminus V}P)^n\un_V[f]\\
 &=&F_f\eq
 where we used that $\un_V=0$ on $M\setminus V$ in the last but one equality.\par
Let $\nu$ be the normalisation into a probability measure of the restriction of $\mu$ to $V$.
\begin{lem}
The probability measure $\nu$ is invariant for $K$.
\end{lem}
\proof
Indeed, we compute that for any $f\in\cB(V)$,
\bq
\nu [K[f]]&=&\frac1{\mu(V)}\mu[\un_V K[f]]\\
&=&\frac1{\mu(V)}\lt(\mu[K[f]]-\mu[\un_{M\setminus V} K[f]]\rt)\\
\eq
By invariance of $\mu$ with respect to $P$, we have
\bq
\mu[\un_{M\setminus V} K[f]]&=&\mu[P[\un_{M\setminus V} K[f]]]\\
&=&\mu\lt[P\un_{M\setminus V}\lt(\sum_{n\in\ZZ_+} (P\un_{M\setminus V})^nP[\un_Vf]\rt)\rt]\\
&=&\mu\lt[\sum_{n\in\NN} (P\un_{M\setminus V})^nP[\un_Vf]\rt]\\
&=&\mu[K[f]]-\mu[P[\un_Vf]]\\
&=&\mu[K[f]]-\mu[\un_Vf]
\eq
In conjunction with the previous identity, we get
\bq
\nu [K[f]]\ =\ \frac1{\mu(V)}\mu[\un_Vf]\ =\ \nu[f]\eq
as wanted.
\wwtbp
\par
From now on, we will only be concerned with the more specific  \textbf{reversible situation} where $P$ is symmetric in $\LL^2(\mu)$
(or equivalently $\mu(dx)P(x, dy) = \mu(dy)P(y, dx)$).
It follows that $P$ can be extended into a bounded self-adjoint operator on $\LL^2(\mu)$.
Then $\nu$ is also reversible with respect to $K$: for any $f,g\in \cB(V)$,
we have
\bq
\nu[fK[g]]&=&\frac1{\mu(V)}\mu[\un_{V}fK[g]]\\
&=&\frac1{\mu(V)}\mu\lt[\un_{V}f \lt(\sum_{n\in\ZZ_+} (P\un_{M\setminus V})^nP[\un_Vg]\rt)\rt]\\
&=&\frac1{\mu(V)}\mu\lt[ \un_Vg \lt(\sum_{n\in\ZZ_+} P(\un_{M\setminus V}P)^n[\un_{V}f]\rt)\rt]\\
&=&\frac1{\mu(V)}\mu\lt[\un_{V}g \lt(\sum_{n\in\ZZ_+} (P\un_{M\setminus V})^nP[\un_Vf]\rt)\rt]\\
&=&\nu[gK[f]]\eq
\par
As a consequence, $K$ can also be extended into a bounded self-adjoint operator on $\LL^2(\nu)$.
It leads us to introduce the following quantities for $k\in\NN$,
\bqn{sigmak}
\sigma_k&\df&\inf_{H\in\cH_k(V)}\sup_{f\in H\setminus\{0\}}\frac{\nu[f(I-K)[f]]}{\nu[f^2]}\eqn
where $\cH_k(V)$ is the set of subspaces of dimension $k$ of $\LL^2(\nu)$. In the above definition and subsequently, the convention $\inf\emptyset\df+\iy$ is enforced.
When $K$ has no essential spectrum, the finite elements of $(\sigma_k)_{k\in\NN}$ are the eigenvalues of $I-K=-S$ with multiplicities,
due to their variational characterization. We want to estimate them via higher order Cheeger inequalities.
To go in this direction, let us consider
\bq
\cA(V)&\df& \{A\in \cA\St A\cap V\in\cA\}\eq
and for $A\in\cA(V)$, the \textbf{Dirichlet--Steklov Markov kernel} $K_A$ defined on $\cB(A\cap V)$ as follows.
For any $f\in \cB(A\cap V)$, consider 
\bq
\fo x\in M,\qquad F_{A,f}(x)&\df& \EE_x[f(Z(\tau_{A\cap V}))\un_{\{\tau_{A\cap V}\leq \tau_{M\setminus A}\}}]\eq
where  $\tau_{A\cap V}$ is the hitting time of $A\cap V$ by $Z$ according to \eqref{tauA}.
The operator $K_A$ is then given by 
\bq
\fo x\in A\cap V,\qquad K_A[f](x)&\df& P[F_{A,f}](x)\eq
\par
Let $\nu_{A}$ be the renormalisation into a probability measure of the restriction of $\mu$ (or $\nu$) to $A\cap V$.
It can be easily checked as above that $K_A$ is Markovian and symmetric in  $\LL^2(\nu_{A})$, so that $K_A$ can be extended into 
bounded self-adjoint operator on $\LL^2(\nu_{A})$.
As in \eqref{sigmak}, we could introduce the quantities  $(\sigma_k(A))_{k\in\NN}$, but only its first element will be important for us:
\bqn{sigmaA}
\sigma_1(A)&\df&\inf_{f\in \LL^2(\nu_{A})\setminus\{0\}}\frac{\nu_A[f(I-K_A)[f]]}{\nu_A[f^2]}\eqn
More precisely, for any $k\in\NN$, let $\cA_k(V)$ be the set of $k$-tuples $(A_1, A_2, ..., A_k)$ of disjoint elements from $\cA(V)$.
We introduce the Dirichlet--Steklov connectivity spectrum $(\kappa_k)_{k\in\NN}$ of $K$ via
\bq
\fo k\in\NN,\qquad 
\kappa_k&\df&\inf_{(A_1, ..., A_k)\in\cA_k(V)}\max_{l\in\lin k\rin} \sigma_1(A_l)\eq
Definition \eqref{sigmaA} can be considered for any $A\in\cA$, but
with the usual convention, we get $\sigma_1(A)=+\iy$ when $A\not\in \cA(V)$, because $\LL^2(\nu_A)=\{0\}$ in this case (and we are left with the trivial $K_A=0$).
Nevertheless, it enables to write 
\bqn{kappak2}
\fo k\in\NN,\qquad 
\kappa_k&=&\inf_{(A_1, ..., A_k)\in\cA_k}\max_{l\in\lin k\rin} \sigma_1(A_l)\eqn
where $\cA_k$ be the set of $k$-tuples $(A_1, A_2, ..., A_k)$ of disjoint elements from $\cA$.
\par\sm
The goal of this section is to show that the extension of Theorem \ref{theo1} holds in this setting:
\begin{theo}\label{theo1m}
There exists a universal constant $c>0$ such that 
\bq
\fo k\in\NN,\qquad  \frac{c}{k^6} \kappa_k\le \sigma_k\,\le\kappa_k \eq
\end{theo}
\par
\sm
As in the finite setting,  the above result leads to higher order Cheeger inequalities presented below.
Nevertheless Theorem \ref{theo1m} is more robust than the latter inequalities \eqref{latterinq1} and \eqref{latterinq2}, as it will appear in its proof.
In a future work, we hope to take advantage of Theorem \ref{theo1m} to give an alternative proof, as well as extensions, of Theorem \ref{C} of the introduction.
\par
We need  the natural extensions of the definitions given in the finite case to our present mesurable state space setting.
The boundary of any $A\in\cA$ is given by
\bq
\pa A&\df& \{(x,y)\St x\in A, \, y\in M\setminus A\}\eq
It is a measurable subset of $M\times M$ endowed with its product $\sigma$-field $\cM\otimes\cM$. 
Consider the measure $\umu$ on $M\times M$  defined by
\bqn{bdry2}
\umu(dx,dy)&=& 
\mu(dx)P(x,dy)\eqn
Here there is a slight difference with the finite case, as we do not impose that the diagonal $D\df\{(x,x)\St x\in M\}$ is negligible with respect to $\umu$:
we cannot do so, because we are not sure $D$ belongs to $\cM\otimes\cM$. It is not important, since we will only integrate with respect to $\umu$ functions
which vanish on the diagonal.
In particular $\umu$ enables to measure $\pa A$ through $\umu(\pa A)$. As a consequence, we can define for $A\in\cA$ the isoperimetric ratios
\bq
\eta(A)&\df& \frac{\umu(\pa A)}{\mu( A)}\\
\eta'(A)&\df& \frac{\umu(\pa A)}{\mu(  A\cap V)}
\eq
(by convention, $\eta'(A)=+\iy$ if $A\not\in\cA(V)$).
Again, the ratio $\eta'(A)$ is the measurable analogue of quantities introduced by Escobar \cite{MR1479552} and Jammes \cite{MR3449183},
since in their terminology, $\pa A$ and $A\cap V$ can be seen respectively as the interior and  exterior boundaries,
when the set $V$ itself is seen as a boundary of $M$.
\par
Next consider
\[\rho(A):=\inf_{\substack{B\in\A\\B\subseteq A}} \eta(B)\]\[\rho'(A):=\inf_{\substack{B'\in\A\\B'\subseteq A} }\eta'(B')\]
\par
For any $k\in\NN$, introduce the $k$-th Cheeger--Steklov constant of $V$ by
\bq
\iota_k&\df&\inf_{(A_1, ..., A_k)\in\cA_k}\max_{l\in\lin k\rin} \rho(A_l)\rho'(A_l)\eq
The next result can be seen as an extension to higher order Cheeger inequalities  of Théorème 1 of Jammes \cite{MR3449183}, as in Theorem {\ref{theo2}:
\begin{theo}\label{theo2m}
Let $c$ be the constant of Theorem \ref{theo1m}.
We have 
\bqn{latterinq1}
\fo k\in\N,\qquad \sigma_k&\geq & \frac{c}{k^6} \iota_k\eqn
\end{theo}
\proof
The deduction of Theorem \ref{theo2m} from Theorem \ref{theo1m} is very similar to that of Theorem \ref{theo2} from Theorem~\ref{theo1}.  For any  function $f\in \LL^2(\nu_A)\setminus\{0\}$, due to Remark \ref{rem10} for the measurable situation  and Lemma \ref{crucial} below,  we have
\bq\frac{\nu_A[f(I-K_A)[f]]}{\nu_A[f^2]}&=&\frac{\mu[F_{A,f}(I-P)[F_{A,f}]]}{\mu[\un_{V\cap A}f^2]}\\
&=&\frac{\int_{M\times M}\mu(dx)P(x,dy)\un_{F_{A,f}(y)\neq F_{A,f}(x)}(F_{A,f}(y)-F_{A,f}(x))^2}{2\mu[\un_{V\cap A}f^2]}
\eq
We multiply the numerator and the denominator by $\int_{M\times M}\mu(dx)P(x,dy)\un_{F_{A,f}(y)\neq F_{A,f}(x)}(F_{A,f}(y)+F_{A,f}(x))^2$ and follow the same calculation as in the proof of Theorem \ref{theo2m}.
The key point is that the statement of the co-area formula is the same in the finite and measurable situations, replacing sums by integrals.
To illustrate the kind of slight modifications  to be taken into account (also  that $\lVe L\rVe$  of Theorem  \ref{theo1} can be replaced by 1 here), let us present the equivalent of the computation \eqref{linkfm}
\bq
\lefteqn{\int_{M\times M}\mu(dx)P(x,dy)\un_{F_{A,f}(y)\neq F_{A,f}(x)}(F_{A,f}(y)+F_{A,f}(x))^2}\\&\leq &2\int_{M\times M}\mu(dx)P(x,dy)\un_{F_{A,f}(y)\neq F_{A,f}(x)}(F_{A,f}^2(y)+F_{A,f}^2(x))\\
&=&4\int_{M\times M}\mu(dx)P(x,dy)\un_{F_{A,f}(y)\neq F_{A,f}(x)}F_{A,f}^2(x)\\
&\leq &4\int_M \mu(dx) F_{A,f}^2(x)\eq
The measurable indicator $\un_{F_{A,f}(y')\neq F_{A,f}(x)}$ is inherited from the Cauchy-Schwarz' inequality in \eqref{linkfm2} and must be kept 
to avoid the possible drawback that $D\not\in\cM\otimes\cM$.
\par
In the same spirit, Definition \eqref{linkfm3} must be replaced by \eqref{linkfm4}. Then we apply the above calculation to a family of functions $f_n\in \LL^2(\nu_A)$ such that $\frac{\nu_A[f_n(I-K_A)[f_n]]}{\nu_A[f_n^2]}\to\sigma_1(A)$ as $n$ tends to $\infty$.
\wwtbp
As in the previous section we consider 
\bqn{latterinq2}
h_k'&\df&\inf_{(A_1, ..., A_k)\in\cA_k(V)}\max_{l\in\lin k\rin} \eta'(A_l)\eqn
and by the same proof,  Proposition \ref{hkup} valid  in the measurable situation, i.e. 
\[\fo k\in\N,\qquad \sigma_k\le  h'_k\]

The proof of Theorem \ref{theo1m} follows the same pattern as in the finite case: it will be deduced from the higher order Cheeger inequalities from \cite{MR3323580},
once the above quantities will be shown to be  limits of spectra associated to speed-up Markov processes.
More precisely, for $r>0$, consider the jump Markov generator $L^{(r)}$ on $M$ given by the kernel
\bq
L\rr(x,dy)&\df& \lt\{
\begin{array}{ll}
r(P(x,dy)-\delta_x(dy))\,,&\hbox{if $x\in M\setminus V$}\\
P(x,dy)-\delta_x(dy)\,,&\hbox{ if $x\in V$}
\end{array}
\rt.\eq
Define
the probability measure $\mu\rr$ on $(M,\cM)$  by
\bq
 \mu\rr(dx)&=&
\lt(\frac{ \un_{V}(x)}{Z_r}+\frac{ \un_{M\setminus V}(x)}{rZ_r}\rt)\mu(dx)
\eq
where $Z_r\df \mu(V)+(1-\mu(V))/r$ is the normalisation constant.\par
The proof of Lemma \ref{lem1} is still valid and leads to
\begin{lem}
The operator $L\rr$ is self-adjoint in $\LL^2(\mu\rr)$.
\end{lem}
\par
Similarly to \eqref{sigmak} and \eqref{sigmaA}, consider 
\bq
\lambda\rr_k&\df&\inf_{H\in\cH_k}\sup_{F\in H\setminus\{0\}}\frac{\mu\rr[F(-L\rr)[F]]}{\mu\rr[F^2]}\eq
where $\cH_k$ is the set of subspaces of dimension $k$ of $\LL^2(\mu)=\LL^2(\mu\rr)$, and for any $A\in\cA$,
\bq
\lambda_1\rr(A)&\df&\inf_{F\in \LL^2(A,\mu)\setminus\{0\}}\frac{\mu\rr[F(-L\rr)[F]]}{\mu\rr[F^2]}\eq
where $\LL^2(A,\mu)$ is the space of $F\in\LL^2(\mu)$ which vanish on $M\setminus A$. The larger $\lambda_1\rr(A)$ is, the easier it is for a (continuous time) Markov process associated to the generator $L\rr$ to exit $A$:
the quantity $\lambda_1\rr(A)$ corresponds to the first Dirichlet eigenvalue of $A$ and measures the asymptotical rate of exit from $A$.
\par
The numerators in the above r.h.s.\ are only slightly dependent on $r\geq 1$ and related to the similar quantities relative to $K$:
\begin{lem}\label{crucial}
We have for any $r>0$ and $F\in\LL^2(\mu)$,
\bq
\mu\rr[F(-L\rr)[F]]&=&\frac1{2Z_r}\int \mu(dx) P(x,dy)(F(y)-F(x))^2\\
&=&\frac1{Z_r}\mu[F(I-P)[F]]
\eq
Furthermore, for any $f\in\LL^2(\nu)$,
\bq
\nu[f(I-K)[f]]&=&\frac1{\mu(V)}\inf\{ \mu[F(I-P)[F]]\St  F_{\vert V}=f\}\\
&=&\frac1{\mu(V)}\mu[F_f(I-P)[F_f]]\eq
where $F_{\vert V}$ stands for the restriction of $F$ to $V$.
\end{lem}
\proof
By definition, for any $r>0$ and $F\in\LL^2(\mu)$, we have
\bq
\mu\rr[f(-L\rr)[F]]&=&-\int_{M\times M} \mu\rr(dx) L\rr(x,dy) F(x)F(y)\\
&=&-\int_{V\times M} \mu\rr(dx) L\rr(x,dy) F(x)F(y)-\int_{(M\setminus V)\times M} \mu\rr(dx) L\rr(x,dy) F(x)F(y)\\
&=&\frac1{Z_r}\int_{V\times M} \mu(dx) (\delta_x(dy)-P(x,dy)) F(x)F(y)\\&&+\frac1{Z_r}\int_{(M\setminus V)\times M} \mu(dx) (\delta_x(dy)-P(x,dy)) F(x)F(y)\\
&=&\frac1{Z_r}\int_{M\times M} \mu(dx) P(x,dy)(F(x)-F(y))F(x)\\
&=&\frac1{2Z_r}\int_{M\times M} \mu(dx) P(x,dy)(F(y)-F(x))^2\eq
where we used the reversibility (under the form $\mu(dx)P(x,dy)=\mu(dy)P(y,dx)$) in the last equality.
Note that the last but one r.h.s.\ is just $\mu[F(I-P)[F]]/{Z_r}$.
\par
Similarly, we compute that for any $f\in\LL^2(\nu)$,
\bq    
\nu[f(I-K)[f]]&=&\int_{V\times V} \nu(dx)K(x,dy) (f(x)-f(y))f(x)\\
&=&\int_V \nu(dx) K[f(x)-f](x)f(x)\\
&=&\int_V \nu(dx) P[f(x)-F_f](x)f(x)\\
&=&\int_{V\times M} \nu(dx)P(x,dy) (f(x)-F_f(y))f(x)\\
&=&\int_{V\times M} \nu(dx)P(x,dy) (F_f(x)-F_f(y))F_f(x)\\
&=&\int_{M\times M} \nu(dx)P(x,dy) (F_f(x)-F_f(y))F_f(x)\\
&=&\frac1{\mu(V)}\mu[F_f(I-P)[F_f]]\eq
where in the last but one equality, we used that $F_f$ is harmonic on $M\setminus V$ according to \eqref{harm}.
It remains to see that
\bqn{harmonicmin}\inf\{ \mu[F(I-P)[F]]\St  F_{\vert V}=f\}&=&\mu[F_f(I-P)[F_f]]\eqn
namely that
among all $F\in\LL^2(\mu)$ coinciding with $f$ on $V$, the quantity
$\mu[F(I-P)[F]]$
 is minimum when $F=F_f$. This is a well-known fact, due to the harmonic property of $F_f$, let us recall the argument.
Write any such function $F$ under the form $F_f+G$ where $G\in\LL^2(\mu)$ vanishes on $V$.
We have
\bq
\mu[F(I-P)[F]]&=&\mu[F_f(I-P)[F_f]]+ \mu[F_f(I-P)[G]]+\mu[G(I-P)[F_f]]+\mu[G(I-P)[G]]
\\
&=&\mu[F_f(I-P)[F_f]]+ 2\mu[G(I-P)[F_f]]+\mu[G(I-P)[G]]\\
&=&\mu[F_f(I-P)[F_f]]+ \mu[G(I-P)[G]]\eq
where we used reversibility, $G_{\vert V}=0$ and \eqref{harm}.
The announced minimisation comes from the non-negativity of
\bq
\mu[G(I-P)[G]]&=&\int_{M\times M} \nu(dx)P(x,dy) (G(x)-G(y))^2\eq
\wwtbp
\par
 Our first approximation results are:
\begin{theo}\label{sgconv}
Assume that $\lambda\df \lambda^{(1)}_1(M\setminus V)>0$ (this quantity will be subsequently called the \textbf{\mbox{Dirichlet} gap} of $M\setminus V$), namely that it is quite easy for the Markov chains $(Z)_{x\in M}$ to enter into $V$.
Then for any $k\in\NN$, we have
\bq
\lim_{r\ri+\iy}\lambda_k\rr&=&\sigma_k\eq
and for any $A\in\cA$,
\bqn{lambda1sigma1}
\lim_{r\ri+\iy} \lambda\rr_1(A)&=&\sigma_1(A)\eqn
More precisely,  the latter convergence is uniform, in the following sense:
let $\frd$ be a distance on the compact set $[0,+\iy]$ compatible with its usual topology.
We have
\bq
\lim_{r\ri+\iy} \sup_{A\in\cA}\frd(\lambda_1\rr(A),\sigma_1(A))&=&0\eq
\end{theo}
More generally, the proof of \eqref{lambda1sigma1} will show that $\lim_{r\ri+\iy} \lambda\rr_k(A)=\sigma_k(A)$,
for any $k\in\NN$, but it will not be useful for our purposes.
\proof
The proof is mainly  concerned with
 the first convergence, since the second one will follow by recycling the obtained quantitative bounds.
\par
We begin by checking that for any $k\in\NN$, we have
\bqn{limsup}
\limsup_{r\ri+\iy}\lambda_k\rr&\leq&\sigma_k\eqn
This result does not require that $\lambda^{(1)}_1(M\setminus V)>0$.
Note that any $H\in\cH_k(V)$ can be seen as an element of $\cH_k$, through the one-to-one mapping
\bq
\LL^2(\nu)\ni f&\mapsto & F_f \in\LL^2(\mu)\eq
so that we have
\bq
\lambda\rr_k&\leq&\inf_{H\in\cH_k(V)}\max_{f\in H}\frac{\mu\rr[F_f(-L\rr)[F_f]]}{\mu\rr[F_f^2]}\eq
According to Lemma \eqref{crucial}, for any $f\in\LL^2(\nu)$,
\bq
\mu\rr[F_f(-L\rr)[F_f]]&=&\frac1{Z_r}\mu[F_f(I-P)[F_f]]\\
&=&\frac{\mu(V)}{Z_r}\nu[f(I-K)[f]]
\eq
Furthermore, we compute that
\bq
\mu\rr[F_f^2]&=&\frac1{Z_r}\lt(\mu[\un_V f^2]+\mu[\un_{M\setminus V} F_f^2]/r\rt)\\
&\geq & \frac1{Z_r}\mu[\un_V f^2]\\
&=& \frac{\mu(V)}{Z_r}\nu[ f^2]
\eq
Thus we get that
\bq
\lambda\rr_k&\leq&\inf_{H\in\cH_k(V)}\max_{f\in H}\frac{\nu[f(I-K)[f]]}{\nu[f^2]}\\
&=&\sigma_k\eq
from which \eqref{limsup} follows at once.\par\sm
Conversely,
to any  subspace $H\subset \LL^2(\mu)$ associate $\wi H$ the subspace of $\LL^2(\nu)$ generated by the functions $F_{\vert V}$ for $F\in H$.
For $k\in\NN$, let $\cH_k^*$ stands for the set of $H\in \cH_k$ which are such that $\wi H\in\cH_k(V)$, namely such that $\wi H$ has dimension $k$.
We begin by remarking that for $k\in\NN$ such that $k\leq \dim(\LL^2(\nu))$ ($\leq+\iy$) and for any $r>0$,
\bqn{Hstar}
\lambda\rr_k&=&\inf_{H\in\cH_k^*}\max_{F\in H\setminus\{0\}}\frac{\mu\rr[F(-L\rr)[F]]}{\mu\rr[F^2]}\eqn
Indeed,
fix some $H\in \cH_k$ and choose $F_1, F_2, ..., F_k$ a basis of $H$. Consider for $l\in\lin k\rin$, $f_l$ the restriction of $F_l$ to $V$.
If $(f_l)_{l\in\lin k\rin}$ is not an independent family of $\LL^2(V)$, then we can find another family $(\wit f_l)_{l\in\lin k\rin}$ of $\LL^2(V)$
such that for any $\epsilon\in(0,1]$, the family $(f_l+\epsilon \wit f_l)_{l\in\lin k\rin}$ is independent.
For  $\epsilon\in(0,1]$, consider $H_\epsilon$ the space generated by $(F_l+\epsilon \wit F_l)_{l\in\lin k\rin}$, where the $\wit F_l$, $l\in\lin k\rin$,
are  the functions coinciding with $\wit f_l$ on $V$ and e.g.\ vanishing outside.
Since $\wi H_\epsilon$ belongs to $\cH_k(V)$, we have that $H_\epsilon\in \cH_k^*$.
Furthermore, it is clear that
\bq
\lim_{\epsilon\ri 0_+} \max_{F\in H_\epsilon\setminus\{0\}}\frac{\mu\rr[F(-L\rr)[F]]}{\mu\rr[F^2]}&=&\max_{F\in H\setminus\{0\}}\frac{\mu\rr[F(-L\rr)[F]]}{\mu\rr[F^2]}
\eq
showing \eqref{Hstar}.\par
Recall that we have by definition
\bq
\lambda&\df& \inf_{F\in\LL^2(M\setminus V,\mu)\setminus\{0\}}\frac{\mu[F(P-I)[F]]}{\mu[F^2]}\\
&=&\inf_{\substack{F\in\LL^2(\mu)\\\un_{M\setminus V}F\neq 0}}\frac{\mu[\un_{M\setminus V}F(I-P)[\un_{M\setminus V}F]]}{\mu[\un_{M\setminus V}F^2]}\eq
It follows that for any $F\in\LL^2(\mu)$, 
\bq
\mu[\un_{M\setminus V}F^2]&\leq &\frac1\lambda \mu[\un_{M\setminus V}F(I-P)[\un_{M\setminus V}F]]\\
&\leq & \frac1\lambda \mu[(F-\un_{V}F)(I-P)[F-\un_{V}F]]\\
&\leq & \frac2\lambda \lt(\mu[F(I-P)[F]]+\mu[\un_{V}F(I-P)[\un_{V}F]]\rt)\\
&\leq & \frac2\lambda \lt(\mu[F(I-P)[F]]+2\mu[\un_{V}F^2]\rt)\eq
where we used that the mapping $\LL^2(\mu)\ni F\mapsto \mu[F(I-P)[F]]$ is a (non-negative) quadratic form (called the \textbf{Dirichlet form} associated to the Markov generator $P-I$, see Remark \ref{rem10})
and that the  spectrum of the  operator $I-P$ is included into $[0,2]$.
We deduce that for any $r>0$,
\bq
\mu\rr[F^2]&=&\frac1{Z_r}\lt(\mu[\un_V F^2]+\frac1r\mu[\un_{M\setminus V} F^2]\rt)\\
&\leq & \frac1{Z_r}\lt(\lt(1+\frac4{r\lambda}\rt)\mu[\un_V F^2]+\frac2{r\lambda}\mu[F(I-P)[F]]\rt)\eq
It follows that
\bq
\frac{\mu\rr[F(-L\rr)[F]]}{\mu\rr[F^2]}&\geq & \frac{\mu[F(I-P)[F]]}{\lt(1+\frac4{\lambda r}\rt)\mu[\un_V F^2]+\frac2{r\lambda}\mu[F(I-P)[F]]}\\
&=& \phi_r\lt(\frac{\mu[F(I-P)[F]]}{\mu[\un_V F^2]}\rt)
\eq
where 
\bq
\phi_r\St [0,+\iy]\ni u&\mapsto& \frac{u}{1+\frac4{\lambda r}+\frac{2u}{\lambda r}}\eq
Note that the latter mapping is increasing, so taking into account Lemma \ref{crucial}, we have, with $f\df F_{\vert V}$,
\bq
 \phi_r\lt(\frac{\mu[F(I-P)[F]]}{\mu[\un_V F^2]}\rt)&\geq &  \phi_r\lt(\frac{\mu[F_f(I-P)[F_f]]}{\mu(V)\nu[f^2]}\rt)
 \\
 &=&\phi_r\lt(\frac{\nu[f(I-K)[f]]}{\nu[f^2]}\rt)
 \eq
 \par
 We deduce from the above computations that for $H\in\cH_k^*$,
 \bq
 \max_{F\in H\setminus\{0\}}\frac{\mu\rr[F(-L\rr)[F]]}{\mu\rr[F^2]}
 &\geq & \max_{f\in \wi H\setminus\{0\}} \phi_r\lt(\frac{\nu[f(I-K)[f]]}{\nu[f^2]}\rt)\\
 &= & \phi_r\lt( \max_{f\in \wi H\setminus\{0\}}\frac{\nu[f(I-K)[f]]}{\nu[f^2]}\rt)\\
&\geq & \phi_r\lt(\sigma_k\rt)
 \eq
since $\wi H\in\cH_k(V)$.
\par
When $k\leq \dim(\LL^2(\nu))$, it follows from \eqref{Hstar} that
\bq
\lambda_k\rr&\geq & \phi_r\lt(\sigma_k\rt)\eq
and it remains to let $r$ go to $+\iy$ to get
\bqn{liminf}
\liminf_{r\ri+\iy}\lambda_k\rr\ \geq\ \lim_{r\ri+\iy}\phi_r(\sigma_k)\ =\ \sigma_k\eqn
\par
When $k>\dim(\LL^2(\nu))$, for any $H\in\cH_k$, we can find $F^*\in H\setminus\{0\}$ such that $F^*_{\vert V}=0$ and so
\bq
 \max_{F\in H\setminus\{0\}}\frac{\mu\rr[F(-L\rr)[F]]}{\mu\rr[F^2]}
 &\geq &\frac{\mu\rr[F^*(-L\rr)[F^*]]}{\mu\rr[F^{*2}]}\\
 &\geq & \phi_r(+\iy)\\
 &=&\lambda r
\eq
It follows that $\lambda_k\rr\geq \lambda r/2$
and letting $r$ go to $+\iy$ we get
\bq
\liminf_{r\ri+\iy}\lambda_k\rr\ =\ +\iy\ =\ \sigma_k\eq
Thus \eqref{liminf} is always true
and in conjunction with \eqref{limsup}, we obtain
the first announced convergence.\par\sm
For the second convergence, note that for $A\in\cA$, the definition of $\sigma_1(A)$ is similar to that of $\sigma_1$
where $V$ is replaced by $V\cup (M\setminus A)$, except we only consider functions that vanish on $M\setminus A$.
It leads us to consider 
\bq
\lambda_A&\df&\lambda_1\uu(A\setminus V)\eq
and for $r>0$, the mapping $\phi_{A,r}$ given by 
\bq
\phi_{A,r}\St [0,+\iy]\ni u&\mapsto& \frac{u}{1+\frac4{\lambda_A r}+\frac{2u}{\lambda_A r}}\eq
The above computations show that for any $r>0$,
\bq
\sigma_1(A)\ \geq \ \lambda_1\rr(A)\ \geq \ \phi_{A,r}(\sigma_1(A))\eq
Note that the mapping $\cA\ni B\mapsto \lambda_1\uu(B)$ is non-increasing with respect to the inclusion of sets
(because $\lambda_1\uu(B)$ corresponds to an infimum over the space of functions $\LL^2(B,\mu)\setminus\{0\}$, which is non-decreasing with respect to $B$),
so we deduce
\bq
\lambda_A&\geq & \lambda\\
\fo r>0,\qquad \phi_{A,r}&\geq &\phi_r
\eq
It follows that to get the wanted uniform convergence, it is sufficient to show that
\bq
\lim_{r\ri+\iy}\sup_{u\in[0,+\iy]}\frd(u,\phi_r(u))&=&0\eq
which is an elementary computation, since it can be reduced to
\bq
\lim_{r\ri+\iy}\max\lt(\sup_{u\in[0,1]}\vert u-\phi(u)\vert, \sup_{u\in[1,+\iy]}\lve \frac1u-\frac1{\phi_r(u)}\rve\rt)&=&0\eq
\wwtbp
\par
\begin{rem}
The assumption of positive Dirichlet gap in Theorem \ref{sgconv} is really needed. 
Indeed, remark that when $\lambda_1^{(1)}(M\setminus V)=0$, then for any $r>0$, we have $\lambda_1\rr(M\setminus V)=0$.
Due to Lemma \ref{crucial}, this is an immediate consequence of
\bq
\fo F\in\LL^2(\mu),\qquad
\frac1{\max(1,r)}\frac{\mu[F(I-P)[F]]}{\mu[F^2]}\ \leq \ 
\frac{\mu\rr[F(-L\rr)[F]]}{\mu\rr[F^2]}\ \leq \ 
\frac1{\min(1,r)}\frac{\mu[F(I-P)[F]]}{\mu[F^2]}\eq
Furthermore, the fact that $\lambda_1\rr(M\setminus V)=0$ implies that $\lambda_2\rr=0$: consider a sequence of functions $(F_n)_{n\in\NN}$ from $\LL^2(M\setminus V,\mu)\setminus\{0\}$ such that
\bq
\lim_{n\ri\iy} \frac{\mu\rr[F_n(-L\rr)[F_n]]}{\mu\rr[F_n^2]}&=&0\eq
and consider for $n\in\NN$, $H_n\df\Vect(\un, F_n)\in\cH_2$. 
We easily get that
\bq\lim_{n\ri\iy} \max_{F\in H_n\setminus\{0\}}\frac{\mu\rr[F(-L\rr)[F]]}{\mu\rr[F^2]}&=&0\eq
i.e.\ $\lambda_2\rr=0$.
In particular, we have
\bq
\lim_{r\ri+\iy}\lambda_2\rr&=&0\eq
But it may happen that $\sigma_2>0$.
Consider for instance an ergodic birth and death transition kernel $P$ on $\ZZ_+$:
we take $M=\ZZ_+$ endowed with a probability measure $\mu$ charging all the points.
The reversible transition kernel $P$ is defined via a Metropolis procedure:
\bq
\fo x,y\in\ZZ_+,\qquad
P(x,y)&\df& \lt\{\begin{array}{ll}
\frac12\lt(\frac{\mu(x)}{\mu(y)}\wedge 1\rt)\,,&\hbox{ if $\lve y-x\rve=1$}\\
0\,,&\hbox{if $\lve y-x\rve\geq 2$}\\
1-\sum_{z\in\ZZ_+\setminus\{x\}}P(x,z)\,,&\hbox{ if $x=y$}
\end{array}\rt.
\eq
where $p\wedge q:=\min\{p,q\}$. The definition of $P$ via the above Metropolis procedure implies that it is irreducible with respect to  $\mu$ (see for example \cite[Section 3.1]{BE04}). Recall that by definition, $P$ is ergodic if and only if $$\fo F\in \LL^2(\mu),\qquad P[F] = F \Rightarrow F\in \Vect(\un)$$
Thus, irreducibility implies ergodicity in the above example. As a result, $P$ is also weakly mixing.
Assume that the queues of $\mu$ are sufficiently heavy, in the sense that
\bq
\lim_{x\ri\iy}\frac{\mu(x)}{\mu([x,\iy))}&=&0\eq
An application of discrete Hardy's inequalities (see \cite{MR1710983}, they are given for finite birth and death processes, but are also valid in the denumerable setting)
implies that
 $\lambda_1^{(1)}(\ZZ_+\setminus\{0,1\})=0$.
 Nevertheless considering for instance $V=\{0,1\}$  we get that $\sigma_2>0$, 
as a consequence of $K(0,1)= P(0,1)>0$ and $K(1,0)= P(1,0)>0$.
More generally it can be proven that $\sigma_2>0$ for  any finite subset of $\ZZ_+$ non-empty and not reduced to a singleton.\par\sm
Note that under the weak mixing assumption (or under the ergodicity assumption), $\lambda_2^{(1)}=0$ means that 0 is the lower bound of the essential spectrum, so that
$\lambda_k^{(1)}=0$ for all $1\leq k< \dim(\LL^2(\mu))+1$ and similarly, $\lambda_k^{(r)}=0$ for any $r>0$ and $1\leq k< \dim(\LL^2(\mu))+1$.
\end{rem}
\par
To prove Theorem \ref{theo1m}  without the assumption of a positive Dirichlet gap on $M\setminus V$, we will accelerate the Markov process associated to the generator $P-I$ more strongly on the slow points of $M\setminus V$
(near $\iy$ in the above remark).
More precisely, we look for a measurable function $\varphi \St M\ri [1,+\iy)$, taking the value 1 on $V$, such that by defining for $r> 0$, the jump Markov generator $L\rr$ by
\bqn{Lr2}
 L\rr(x,dy)&\df& \lt\{
\begin{array}{ll}r\varphi(x)(P(x,dy)-\delta_x(dy))\,,&\hbox{ if $x\in M\setminus V$}\\
\varphi(x)(P(x,dy)-\delta_x(dy))\,,&\hbox{ if $x\in V$}
\end{array}\rt.
\eqn
we have that $L^{(1)}$ admits a positive Dirichlet gap on $M\setminus V$.
Then, with the corresponding spectra, Theorem \ref{sgconv} will hold.
Note that the notions of harmonic functions on $M\setminus V$ with respect to $P-I$ and $L\rr$, for all $r>0$,  coincide  and the corresponding Steklov Markov kernels and generators are the same.
\par
Let $X\df(X(t))_{t\geq 0}$ be a jump Markov process of generator $P-I$ (see Chapter 4 in \cite{MR838085} for the definition). Fix some $\chi\in (0,1)$ and consider the function $\varphi$ defined by
\bq
\fo x\in M,\qquad \varphi(x)&\df& \frac{1}{\EE_x[\chi^{\tau}]}\eq
where $\tau\df\inf\{t\geq 0\St X_t\in V\}$.
Note that when $x\in M \setminus V$ is a point from which it is difficult to hit $V$, namely such that $\tau$ has a propensity to be large, then $\varphi(x)$ is quite large also: the jump Markov process 
$X\uu\df(X\uu(t))_{t\geq 0}$ associated to $L\uu$ is strongly accelerated at $x$ in comparison with $X$, as wanted. 
From now on, the notation $L\rr$, for $r>0$, will only refer to the operators given in \eqref{Lr2}.
Here is the consequence of the acceleration procedure:
\begin{lem}\label{Etau1}
We have
\bq
\fo x\in M,\qquad \EE_x[\tau\uu]&\leq &  \frac1{\ln(1/\chi)}\eq
where $\tau\uu\df\inf\{t\geq 0\St X\uu_t\in V\}$
\end{lem}
\proof
Let us recall the time change transformations (cf.\ for instance Chapter 6 from the book of Ethier and Kurtz \cite{MR838085}), which enable to construct $X\uu$ from $X$ when both processes start from a fixed $x\in M$.
Due to \cite[Theorem 1.4]{MR838085}, if we define $(\theta_t)_{t\geq 0}$ via
\bq
\fo t\geq 0,\qquad 
\int_0^{\theta_t}\frac{1}{\varphi(X_s)}\, ds&=&t\eq
then we can take
\bq
\fo t\geq 0,\qquad X\uu(t)&\df& X(\theta_t)\eq
In particular, we get
\bq
\tau\uu&=&\int_0^\tau \frac{1}{\varphi(X_s)}\, ds\eq
It follows that
\bq
\EE_x[\tau\uu]&=&\EE_x\lt[\int_0^\tau \frac{1}{\varphi(X_s)}\, ds\rt]\\
&=&\int_0^{+\iy}\EE_x\lt[\un_{s\leq \tau}\frac{1}{\varphi(X_s)}\rt]\, ds\\
&=&\int_0^{+\iy}\EE_x\lt[\un_{s\leq \tau}\EE_{X_s}[\chi^{\tau}]\rt]\, ds\\
&=&\int_0^{+\iy}\EE_x\lt[\un_{s\leq \tau}\chi^{-s}\EE_{x}[\chi^{\tau}\vert (X_u)_{u\in[0,s]}]\rt]\, ds\\
&=&\int_0^{+\iy}\EE_x\lt[\un_{s\leq \tau}\chi^{-s}\chi^{\tau}\rt]\, ds\eq
where we use the measurability  of the event $\{s\leq \tau\}$ with respect
to the $\sigma$-field generated by $(X_u)_{u\in[0,s]}$, the fact that on $\{s\leq \tau\}$, we have $\tau=s+\tau\circ \theta_s$, where $\theta_s$ is the shift of the trajectories by an amount  $s$ of time,
and the Markov property, stating that for any measurable functional $F$ on the trajectories, we have a.s. $\EE_x[F\circ \theta_s\vert (X_u)_{u\in[0,s]}]=\EE_{X_s}[F]$. In this formula, $\EE_{X_s}$ is the  expectation with respect to a diffusion $X$ starting from $X_s$ at time 0.
Since all the integral elements are non-negative, we can use again Fubini's formula to get
that the last integral is equal to
\bq
\EE_x\lt[\int_0^{+\iy}\un_{s\leq \tau}\chi^{-s}\chi^{\tau}\, ds\rt]&=&\EE_x\lt[\int_0^{\tau}\chi^{\tau-s}\, ds\rt]\\
&=&\EE_x\lt[\int_0^{\tau}\chi^{s}\, ds\rt]\\
&=&\EE_x\lt[\frac{\chi^{\tau}-1}{\ln(\chi)}\rt]\\
&\leq & \frac1{\ln(1/\chi)}\eq
as announced.\wwtbp\par
From the previous uniform boundedness of the expectations of $\tau\uu$, we deduce uniform exponential bounds on its queues:
\begin{lem}\label{explam}
We have
\bq
\fo x\in M,\,\fo s\geq 0,\qquad 
\PP_x[\tau\uu\geq s]&\leq & 2\exp(-\alpha s)\eq
with $\alpha\df \ln(2)\ln(1/\chi)/2$.
\end{lem}
\proof
For any $n\in\ZZ_+$, we have
\bq
\fo x\in M,\qquad \PP_{x}[\tau\uu\geq an ]&\leq & 2^{-n}\eq
where
\bq
a&\df&\frac{2}{\ln(1/\chi)}\eq
This is shown by iteration on $n\in\ZZ_+$. It is clear for $n=0$ and if it is true for some $n\in\ZZ_+$,
then by the Markov property and Lemma \ref{Etau1}: for any $x\in M$,
\bq
\PP_{x}[\tau\uu\geq a(n+1)]&=&\EE_x[\un_{\tau\uu\geq a}\PP_{X\uu(a)}[\tau\uu\geq an]]\\
&\leq & 2^{-n}\PP_x[\tau\uu\geq a]\\
&\leq & 2^{-n}\frac{\EE_x[\tau\uu]}{a}\\
&\leq & 2^{-n}\frac{1}{a\ln(1/\chi)}\\
&=&2^{-(n+1)}\eq
where in the third line we use the Markov inequality.
\par
For any $s\in\RR_+$, write $n\df \lfloor s/a\rfloor\in\ZZ_+$, so that
\bq
\fo x\in M,\qquad
\PP_x[\tau\uu\geq s]&\leq & \PP_x[\tau\uu\geq na]\\
&\leq &2^{-n}\\
&=&2^{-\lfloor s/a\rfloor}\\
&\leq &2(2^{-s/a})\\
&=&2\exp(-\alpha s)\eq
as announced.\wwtbp
\par
To simplify the notation, we now take $\chi=\exp(-2/\ln(2))$, so that $\alpha=1$.
Uniform exponential bounds on the queues of exit times from a domain are well-known to imply that the associated Dirichlet gap is positive.
Here is a simple proof of this fact:
\begin{lem}
We have \bq
\lambda^{(1)}_1(M\setminus V)&\geq &\frac12 \eq
where the l.h.s.\ is relative to the accelerated generator $L^{(1)}$.
\end{lem}
\proof
As in Lemma \ref{lem1}, we see that the measure $\frac1{\varphi(x)}\mu(dx)$ is reversible for $L\uu$.
Its total weight is 
\bq
Z\uu\ \df\ \int \EE_x[\chi^{\tau\uu}]\, \mu(dx)\ \in\ (0,1)\eq
which leads us to define $\mu\uu(dx)\df \frac1{Z\uu\varphi(x)}\mu(dx)$, the invariant probability for $L\uu$.\par
Our goal is to show that
\bqn{D12}
\nonumber\lambda_1\uu(M\setminus V)&\df&\inf_{F\in \LL^2(M\setminus V,\mu\uu)\setminus\{0\}}\frac{\mu\uu[F(-L\uu)[F]]}{\mu\uu[F^2]}\\
&\geq &\frac12  \eqn
\par
So consider $F$ a bounded and measurable function on $M$, vanishing on $V$.
By the martingale problems associated to $X\uu$, there exists a $\LL^2$ martingale $(M_t)_{t\geq 0}$
such that
\bq
\fo t\geq 0,\qquad
F^2(X\uu(t))&=&F^2(X\uu(0))+\int_0^t L\uu[F^2](X\uu(s))\, ds +M_t\eq
Replace in this relation $t$ by $t\wedge \tau\uu$ and take the expectation to get
\bq
\EE[F^2(X\uu(t\wedge \tau\uu))]&=&\EE[F^2(X\uu(0))]+\EE\lt[\int_0^{t\wedge \tau\uu} L\uu[F^2](X\uu(s))\, ds\rt]\eq
where we use the martingale property $\EE(M_t)=\EE(M_0)=0$. Via dominated convergence, we can let $t$ go to infinity to obtain
\bq
\EE[F^2(X\uu( \tau\uu))]&=&\EE[F^2(X\uu(0))]+\EE\lt[\int_0^{\tau\uu} L\uu[F^2](X\uu(s))\, ds\rt]\eq
Note that since $X\uu( \tau\uu)\in V$ the l.h.s.\ vanishes, we deduce
\bq
\EE[F^2(X\uu(0))]&=&-\EE\lt[\int_0^{\tau\uu} L\uu[F^2](X\uu(s))\, ds\rt]\eq
We have not yet specified the initial distribution of $X\uu(0)$, but take it now to be $\mu\uu$,
so the l.h.s.\ is
\bq
\EE_{\mu\uu}[F^2(X\uu( 0))]&=&\int \mu\uu(dx) F^2(x)=\mu\uu[F^2] \\
\eq
Concerning the r.h.s., recall that the \textbf{carré du champs} $\Gamma\uu$ associated to $L\uu$ and defined on any bounded and measurable function $G$ on $M$ by
\bq
\Gamma\uu[G]&\df& L\uu[G^2]-2GL\uu[G]\eq
is a non-negative function (cf.\ for instance the book of Bakry, Gentil and Ledoux \cite{MR3155209}).
It follows that
\bq
-\EE_{\mu\uu}\lt[\int_0^{\tau\uu} L\uu[F^2](X\uu(s))\, ds\rt]&\leq &-2\EE_{\mu\uu}\lt[\int_0^{\tau\uu} F(X\uu(s))L\uu[F](X\uu(s))\, ds\rt]\\
&\leq & 2\EE_{\mu\uu}\lt[\int_0^{\tau\uu}\vert F(X\uu(s))L\uu[F](X\uu(s))\vert\, ds\rt]\\
&=&\int_0^{+\iy}\EE_{\mu\uu}\lt[\un_{s\leq \tau\uu}\vert F(X\uu(s))L\uu[F](X\uu(s))\vert\rt]\, ds\eq
For any $s\geq 0$, taking into account
Lemma \ref{explam}, we have
\bq
\EE_{\mu\uu}\lt[\un_{s\leq \tau\uu}\vert F(X\uu(s))L\uu[F](X\uu(s))\vert\rt]&=&\EE_{\mu\uu}\lt[\PP_{X\uu(s)}[s\leq \tau\uu]\vert F(X\uu(s))L\uu[F](X\uu(s))\vert\rt]\\
&\leq & 2\exp(-s)\EE_{\mu\uu}\lt[\vert F(X\uu(s))L\uu[F](X\uu(s))\vert\rt]\\
&= & 2\exp(-s)\mu\uu[\vert F L\uu[F]\vert]\eq
where we used the invariance of $\mu\uu$ (meaning that for any $s\geq 0$,  the law of $X\uu(s)$ is equal to $\mu\uu$ when the initial law is $\mu\uu$).
We have thus proven that
\bq
\mu\uu[F^2]&\leq & \int_0^{+\iy} 2\exp(-s)\mu\uu[\vert F L\uu[F]\vert]\, ds\\
&=&2\mu\uu[\vert F L\uu[F]\vert]\\
&\leq & 2\sqrt{\mu\uu[F^2]\mu\uu[(L\uu[F])^2]}
\eq
i.e.
\bq
\mu\uu[F^2]&\leq & 4\mu\uu[(L\uu[F])^2]\eq
The fact that $L\uu$ is a non-positive self-adjoint operator enables to see that this relation extend to any function in the domain of $L\uu$ with Dirichlet condition on $V$.
It follows that the spectrum of $-L\uu$ with Dirichlet condition on $V$ is above 1/2, which amounts to \eqref{D12}.\wwtbp
\par
As already mentioned, the Steklov Markov kernel $K\uu$ associated to $L\uu$ and $V$ is the same as $K$.
Since in general the generator $L\uu$ cannot be written under the form $P\uu-I$, where $P\uu$ would be a Markov kernel on $M$,
the definitions \eqref{Fgen} and \eqref{Skern} must be slightly generalized: denote for any $f\in\cB(V)$,
\bqn{Fgen2}
\fo x\in M,\qquad F\uu_f(x)&\df& \EE_x[f(X\uu(\tau\uu))]\\
\nonumber\fo x\in V,\qquad K\uu[f](x)&\df& L\uu[F_f\uu](x)+f(x)\eqn
where $\tau\uu$ was defined in Lemma \ref{Etau1}.
The latter expression for $K\uu$ may appear strange at first view; it is due to the fact that 
 it is a Markov kernel operator. If we rather consider the Steklov generator $S\uu\df K\uu-I$,
we get the more natural formulation: $S\uu[f]=L\uu[F_f]$, for $f\in\cB(V)$, as in the section on finite Markov process.
Coming back to our previous convention of Steklov Markov kernels,
note that for any $x\in V$, we have
\bq
 L\uu[F_f\uu](x)+f(x)&=&L[F_f\uu](x)+F_f\uu(x)\\
 &=&\int F_f\uu (y)\,P(x,dy)\\
 &=&P[F_f\uu](x)\eq
 more in adequacy with \eqref{Skern}.
Note furthermore that the function $F_f\uu$ defined by \eqref{Fgen2} is the $L\uu$-harmonic extension of $f$ to $M$: it satisfies
\bq
\lt\{\begin{array}{ll}
L\uu[F_f]=0\,,&\quad\hbox{on $M\setminus V$}\\
F_f\uu= f\,,&\quad\hbox{ on $V$}
\end{array}\rt.\eq
Since $L\uu=\varphi L$, with $\varphi$ non-vanishing, the condition $L\uu[F_f]=0$ is the same as $L[F_f]=0$.
It follows that $F_f\uu=F_f$ and finally $K\uu[f]=K[f]$.
By completion, this is true on $\LL^2(\nu)$, i.e.\ $K\uu=K$.
The equality $F_f\uu=F_f$ is also obvious from the probabilistic point of view, since $X\uu$ is a time change of $X$ (as seen in the proof of Lemma \ref{Etau1}),
which itself is the Poissonisation of the Markov chain $Z$ with the same initial condition and associated to $P$:
let $(\cE_n)_{n\in\NN}$ be independent exponential random variables of parameter 1, $X$ can be constructed from $Z$ via
\bq
\fo t\geq 0,\qquad X_t&=& Z_n,\quad \hbox{where $n\in\ZZ_+$ is such that
$\sum_{p=1}^{n}\cE_p\leq t<\sum_{p=1}^{n+1}\cE_p$}\eq
\par
The previous considerations are also valid  for the operators $K_A\uu$, defined in a similar fashion for $A\in \cA(V)$ and we get that $K_A\uu=K_A$.
We can now apply Theorem \ref{sgconv} with respect to the generator $L\uu$, which by construction admits a Dirichlet gap on $M\setminus V$.
The l.h.s.\ in the two convergences of Theorem~\ref{sgconv} correspond to the generators given by \eqref{Lr2} and the 
r.h.s.\ are  given by \eqref{sigmak} and \eqref{sigmaA}, according to the above discussion.
These convergences are our final approximation results for  the  quantities $(\sigma_k)_{k\in\NN}$ and $(\sigma_1(A))_{A\in\cA}$ .
\par\sm
We can now come to the
\prooff{Proof of Theorem \ref{theo1m}}
The upper bound is an immediate consequence of the definition of $\sigma_k$. Indeed for every $(A_1, ..., A_k)\in\cA_k$ it is enough to consider the vector space generated by a family $\{f_{l,n}\in L^2(A_l,\mu)\St l\in\lin k\rin\}$ of test functions such that $\frac{\nu_{A_l}[f_{n,l}(I-K_{A_l})[f_{n,l}]]}{\nu_{A_l}[f_{n,l}^2]}$ tends to $\sigma_1(A_l)$ as $n\to\infty$.\\ For the lower bound, similarly to \eqref{kappak2}, define for any $r>0$,
\bq
\fo k\in\NN,\qquad 
\Lambda_k\rr&=&\inf_{(A_1, ..., A_k)\in\cA_k}\max_{l\in\lin k\rin} \lambda\rr_1(A_l)\eq
We have seen in \cite{MR3323580}, extending the similar result  Lee, Oveis Gharan and Trevisan \cite{MR2961569} gave in a finite setting,
that there exists a universal constant $c>0$ such that 
\bqn{LGTM}
\fo r>0,\,\fo k\in\NN,\qquad \lambda\rr_k&\geq & \frac{c}{k^6} \Lambda\rr_k\eqn
Fix some $k\in\NN$.
The first convergence of Theorem \ref{sgconv} shows that the l.h.s.\ converges to $\sigma_k$ as $r$ goes to $+\iy$.
Its uniform convergence leads to
\bq
\lim_{r\ri+\iy} \Lambda_k\rr&=&\kappa_k\eq
so we can pass to the limit in \eqref{LGTM} to obtain the announced inequality.
\wwtbp
We end this section with Proposition \ref{proA} in the introduction. 
\begin{pro}\label{ms-1} There is a universal positive constant $c'$  such that
\bq\fo k\in \N,\qquad \sigma_{2k}&\ge&
\frac{c'}{\log^2(k+1)}\iota_k\eq
\end{pro}
\proof By \cite{MR3323580}, the proof of Proposition \ref{fs-1}  can be extended here. In particular, we have 
\bq\fo k\in \N,\qquad \lambda^{(r)}_{2k}&\ge&
\frac{c}{\log^2(k+1)}\Lambda^{(r)}_k\eq
and
\bq\fo k\in\lin v\rin,\qquad \Lambda_k&\geq  &\frac{1}{8}\iota_k\eq
and the statement follows. 

\wwtbp\par
\section{The  Riemannian manifold framework}\label{RM}
Let $(M,g)$ be a compact Riemannian manifold of dimension $n$ with smooth boundary. We assume that $M$ is connected. Recall  the Steklov problem \eqref{a-stek} considered in the introduction:
\bq
\left\{\begin{array}{ll}
\Delta f=0\,,&\quad\hbox{in}\;\; M\\
\frac{\partial f}{\partial \nu}=\sigma f\,,&\quad\hbox{on} \;\; \partial M
\end{array}\right.\eq
where $\nu$ is the unit outward normal to the boundary. Our goal, as in the previous sections, is to relate its eigenvalues $0=\sigma_1<\sigma_2\le\cdots\le\sigma_k\le\cdots\nearrow\infty$ to some isoperimetric constants.  
  We first show that that \eqref{a-stek} can be seen as a limit of a family of Laplace eigenvalue problems.   This   is already known due to the results of Lamberti and Provenzano \cite{LP15,Pthese}. They showed that  the Steklov eigenvalue problem \eqref{a-stek} can be considered as the limit of the family of Neumman eigenvalue problems 
\begin{equation}\label{a-neu}\left\{\begin{array}{ll}
\Delta f+\lambda\rho_\epsilon f=0\,,&\quad \hbox{in}\;\;  M\\
\frac{\partial f}{\partial \nu}=0\,,&\quad \hbox{on} \;\; \partial M
\end{array}\right.\end{equation}
for $\epsilon$  small enough (one can choose $\epsilon$ for example smaller than the focal distance of $\partial M$).  Here  $ M_\epsilon:=\{x\in M:\; d(x,\partial M)<\epsilon\}$, and 
\begin{equation}\label{rhodef}\rho_\epsilon(x)=\left\{\begin{array}{ll}
\epsilon\,, &\quad\hbox{if}\;\, x\in M\setminus { M_\epsilon}\\
\frac{1}{\epsilon}\,,&\quad\hbox{if}\,\, x\in  M_\epsilon
\end{array}\right.\end{equation}
We denote the eigenvalues of problem \eqref{a-neu} by 
$$0=\lambda^\epsilon_1<\lambda^\epsilon_2\le\cdots\le\lambda^\epsilon_k\le\cdots\nearrow\infty$$
Then we have 
\begin{theo}\cite{LP15,Pthese}\label{provenzano}
For every $k\in\N$
\bqn{stek-neu}\lim_{\epsilon\to0}\lambda_k^\epsilon&=&\sigma_k\eqn
\end{theo}
\begin{rem}  
We remark that Lamberti and Provenzano \cite{LP15,Pthese} stated the above convergence for bounded domains in $R^n$ with smooth boundary, and the definition of $\rho_\epsilon$ on $\partial M$ is slightly different. However, a verbatim proof also results   in the  convergence \eqref{stek-neu} on a compact  Riemannian manifold $(M,g)$ with smooth boundary, see \cite[Chapter 3]{Pthese} for the details of the proof.   \end{rem}
One can see the similarity of the above theorem with the statement  of  Proposition \ref{pro2} and Theorem~\ref{sgconv}. It would be very interesting to have an alternative approach to prove Theorem \ref{provenzano} and Theorem \ref{hm-stek} below by using the results of the previous section. We hope to obtain a unified approach  in a future work.\\

Let  $ A\subset  M$ be a nonempty open domain in $M$. Let $\partial_eA:=\bar A\cap \partial M$ and $\partial_i A:=\partial A\cap \Int\, M$ be smooth manifolds of  dimension $n-1$ when they are nonempty sets.  
We consider the  mixed Dirichlet--Sobolev eigenvalue problem
 \begin{equation}\label{dir-stek}\left\{\begin{array}{ll}
\Delta  f=0&in\;\; A\\
\frac{\partial f}{\partial \nu}=\sigma f&on \;\; \partial_eA\\
 f=0&on \;\; \partial_iA\\
\end{array}\right.\end{equation}
We also need to consider the following mixed Dirichlet--Neumann eigenvalue problem
 \begin{equation}\label{dir-neu}\left\{\begin{array}{ll}
\Delta  f+\lambda\rho_\epsilon f=0&in\;\; A\\
\frac{\partial f}{\partial \nu}=0&on \;\; \partial_eA\\
 f=0&on \;\; \partial_iA\\
\end{array}\right.\end{equation}
where $\rho_\epsilon$ is defined in \eqref{rhodef}.\\
 
If $\partial_i A=\emptyset$, then $A=\Int\,M$ and the first eigenvalue is zero. Otherwise the first eigenvalues of  the eigenvalue problem \eqref{dir-stek} and \eqref{dir-neu} are not zero  and we denote their eigenvalues  by 
$$0<\sigma_1(A)\le\sigma_2(A)\le\cdots\le\sigma_k(A)\le\cdots\nearrow\infty$$ 
and 
$$0<\lambda^\epsilon_1(A)\le\lambda^\epsilon_2(A)\le\cdots\le\lambda^\epsilon_k(A)\le\cdots\nearrow\infty$$ 
respectively. When $\partial_eA=\emptyset$, our  convension is that $\sigma_k(A)=\infty$, for every $k\in \N$.
%
 Denote by $\mathcal{A}$  the set of nonempty open domains in $M$ such that $\partial_iA$ and $\partial_e A$ are smooth  sub-manifolds of dimension $n-1$ when they are nonempty.  Let $\A_k$ be the set of $k$-tuple $(A_1,...,A_k)$ of mutually disjoint elements of $\A$. We define
\bqn{minmax1}\Lambda_k^\epsilon&:=&\inf_{(A_1,...,A_k)\in\A_k}\max_{l\in \llbracket k \rrbracket} \lambda_1^\epsilon(A_l) \eqn
The higher order Cheeger inequality for eigenvalues $\lambda_k^\epsilon(M)$, $k\in \N$ was proved by Miclo in \cite{MR3323580}:
\begin{theo}[\cite{MR3323580}]\label{miclo}There exists a universal constant $c>0$ such that for any compact Riemannian manifold $M$ with smooth boundary, the eigenvalues  $\lambda_k^\epsilon$ of Neumann eigenvalue problem \eqref{a-neu} satisfy
 \[\frac{c}{k^6}\Lambda_k^\epsilon\le\lambda_k^\epsilon\le \Lambda_k^\epsilon\qquad \fo  k\in\N\]
\end{theo} 
\begin{rem}
The above theorem in \cite{MR3323580} is originally stated for the Laplace eigenvalue problem with smooth coefficients on closed manifolds. But the argument remains the same when we consider the Neumann eigenvalue problem \eqref{a-neu} on a compact manifold with smooth boundary. 
\end{rem}
Similar to Defintion \eqref{minmax1}, we define 
\bq\kappa_k&:=&\inf_{(A_1,...,A_k)\in\A_k}\max_{l\in \llbracket k \rrbracket} \sigma_1(A_l)\eq
\begin{theo}\label{hm-stek}
There exists a universal constant $c_1$ such that for any compact Riemannian manifold $M$ with boundary and for any $k\in \N$, the eigenvalues  $\sigma_k( M)$ of problem \eqref{a-stek} satisfy
 \[\frac{c_1}{k^6}\kappa_k\le\sigma_k\le\kappa_k\]
\end{theo}
As a consequence of Theorem \ref{hm-stek} we get the higher order Cheeger--Steklov inequalities, see Theorem \ref{cheeg-stek} below. We first define the Cheeger--Steklov constants in this setting similar to those already discussed in the previous sections.  
For any open subset $A$ of $M$ with piecewise smooth boundary,  let $\mu(A)$ denote its  Riemannian measure and $\umu(\partial A)$ be the induced $(n-1)$-dimensional Riemannian measure of $\partial A$. We define for every $A\in \A$ the isoperimetric ratios 
\bq
\eta(A)&\df& \frac{\umu(\pa_i A)}{\mu( A)}\\
\eta'(A)&\df& \frac{\umu(\pa_iA)}{\umu(   \partial_e A)}
\eq
 Note that $\eta'(A)=\infty$ if $\bar A\cap \partial M=\emptyset$. Let
\bqn{rhoa}\rho(A)&:=&\inf_{\substack{B\in\A\\ B\subset A\\ \bar B\cap \partial_iA=\emptyset}} \eta(B)\eqn
\bq\rho'(A)&:=&\inf_{\substack{B'\in\A\\ B'\subset A\\ \bar B'\cap \partial_iA=\emptyset}} \eta'(B')\eq
 For any $k\in\N$ we define the $k$-th Cheeger--Steklov constant of $M$ by 
\bq\iota_k&:=&\inf_{(A_1,\cdots,A_k)\in\A_k}\max_{l\in\llbracket k\rrbracket}\rho(A_l)\rho'(A_l).\eq

The following theorem extends the results of Escobar \cite{MR3323580} and Jammes \cite{MR3449183}.
\begin{theo}\label{cheeg-stek}
There exists a universal constant $c$ such that for any compact Riemannian manifold $M$ with smooth boundary and for any $k\in \N$, the eigenvalues  $\sigma_k( M)$ of problem \eqref{a-stek} satisfy
 \bq\sigma_k&\ge&\frac{c}{k^6}\iota_k\eq
\end{theo}
\begin{rem}\label{iotaseq}\begin{itemize}
 \item[i)] One can check that for every $k\in \N$ one has $\iota_k\le\iota_{k+1}$. This  is also true in finite and measurable situation.
 \item[ii)] Note that $\eta'(B)$ is scale invariant. Hence, as  mentioned in \cite{MR3449183}, the power of $\eta(B)$  has to be one so that $\iota_k$ has the same scaling as $\sigma_k$. 
 \end{itemize}
\end{rem}

Note that for $k=2$, Theorem \ref{cheeg-stek} gives a version of Jammes' result \cite{MR3449183}. 
The above theorem is the  direct sequence  of Theorem \ref{hm-stek} and Lemma \ref{0stek} below. 
\begin{lem}\label{0stek}Let $\sigma_1(A)$ be the first eigenvalue of the Dirichlet-Steklov eigenvalue problem \eqref{dir-stek}. Then we have
\bq\sigma_1(A)&\ge&\frac{1}{4}\rho(A)\rho'(A)\eq
\end{lem}
\proof
 Let $f$ be the eigenfunction associated with $\sigma_1(A)$. We repeat  the same argument as Jammes' argument in \cite{MR3449183} to estimate $\sigma_1(A)$.
\bq\sigma_1(A)&=&\frac{\int_A|\nabla f|^2\,{d\mu}\,\int_Af^2d\mu  }{\int_{\partial_eA}f^2d\umu \,\int_Af^2d\mu  }\ge\frac{\left(\int_{A}|f\nabla f|d\mu  \right)^2}{\int_{\partial_eA}f^2d\umu \,\int_Af^2d\mu  }\\
&\ge&\frac{1}{4}\left(\frac{\int_{A}|\nabla f^2|d\mu  }{\int_{\partial_eA}f^2d\umu }\right)\left(\frac{\int_{A}|\nabla f^2|d\mu  }{\int_Af^2d\mu  }\right)\eq
where $d\mu  $ and $d\umu $ are $n$-dimensional and $(n-1)$-dimensional Riemannian volume elements respectively. \\
Let $h:=f^2$ and $H_t:=h^{-1}[t,\infty)$. Note that $H_t\in\cA$ almost surely in $t$. Then by the co-area formula we have
\begin{equation*}\left(\frac{\int_{A}|\nabla h|d\mu  }{\int_{\partial_eA}h\,d\umu }\right)\left(\frac{\int_{A}|\nabla h|d\mu  }{\int_Ah\,d\mu  }\right)=\left(\frac{\int_{0}^\infty\umu(\partial_iH_t)dt}{\int_0^\infty \umu(\partial_eH_t)\,dt }\right)\left(\frac{\int_{0}^\infty\umu(\partial_iH_t)dt}{\int_0^\infty\mu(H_t)\,dt  }\right)
\ge\rho(A)\rho'(A)\end{equation*}
which completes the proof.
\wwtbp\par
It  remains to prove Theorem \ref{hm-stek}.
\prooff{Proof of Theorem \ref{hm-stek}}
Recall that by the variational characterisation of Steklov eigenvalues 
\bq\sigma_k&\le&\max_{j\in\llbracket k\rrbracket} \frac{\cE_\Delta({f_j,f_j})}{\int_{\partial M}f_j^2d\umu }\eq
where $\{f_j\}$ is a family of test functions in $H^1(M)$ with mutually disjoint supports and $\cE_\Delta(f,f):=\int_M|\nabla f|^2d\mu $ is the Dirichlet form associated to $\Delta$. Hence,
the upper bound of $
\sigma_k$ is a direct consequence of the variational characterisation of Steklov eigenvalues. 
 
We now prove the lower bound.  We need the following key lemma. 

\begin{lem} \label{keylem}{The following inequality holds.}\bq\lim_{\epsilon\to0}\Lambda_k^\epsilon&\ge&\frac{1}{4} \kappa_k\eq
\end{lem}

\proof
Let $(A_1,\cdots,A_k)\in\mathcal{A}_k$ and $H^1_0(A_j,\partial_i A_j)$ be  the closure of $\{f\in C^\infty(A_j)\,:\, f\equiv 0\;\hbox{on $\partial_iA_j$}\}$ in $H^1(A_j)$. {We can assume $\partial_eA_j\neq\emptyset$.} For any $\epsilon$ small enough (will be determined {below})  and every $f\in H^1_0(A_j,\partial_i A_j)$, $j\in\llbracket k\rrbracket$ we give an upper bound for the denominator of 
\bqn{psi}\frac{\int_{A_j} |\nabla f|^2d\mu  }{\int_{A_j}\rho_\epsilon f^2\,d\mu  }&=&\frac{\int_{A_j} |\nabla f|^2d\mu  }{\frac{1}{\epsilon}\int_{A_{j,e}^\epsilon} f^2\,d\mu  +\epsilon \int_{A_j\setminus A_{j,e}^\epsilon} f^2\,d\mu  }\eqn%
where $A^\epsilon_{j,e}:=\{x\in A_j: d(x,\partial M)<\epsilon\}$. For every $f\in H^1_0(A_j,\partial_i A_j)$  consider $\un_{A_{j}}f$ as an element of $H^1(M)$. Then
\bq\frac{1}{\epsilon}\int_{A_{j,e}^\epsilon} f^2\,d\mu & =&\frac{1}{\epsilon}\int_{M_\epsilon} \un_{A_{j}}f^2\,d\mu   \eq
There exists $\epsilon_0>0$ such that  for every $\epsilon\in(0,\epsilon_0)$ the map
\begin{eqnarray*}E:\partial M\times(0,\epsilon)&\to& M_\epsilon\\
(x,t)&\mapsto& \exp_x(-t\nu(x))
\end{eqnarray*}
 is a diffeomorphism. Note that $|\det DE(x,t)|=1+O(t)$. Hence, by choosing $\epsilon_0$ even smaller, we can impose that for all $(x,t)\in \partial M\times(0,\epsilon)$
\[\quad \sup_{s\in(0,t)}\frac{|\det DE(x,t)|}{|\det DE(x,s)|}\le 2,\quad\mbox{which also implies},\quad \vert\det DE(x,t)\vert\le2\]
Let  $F\in H^1(M) $ and by abuse of notation, denote $F\circ E$ by $F$. For a.e.  $(x,t)\in \partial M\times(0,\epsilon) $ we have 
\bq |F(x,t)|&\le&  |F(x,0)|+\int_0^t\left|\frac{\partial  F}{\partial s}(x,s)\right|ds\eq
Thus
\bq\frac{1}{\epsilon}\int_{M_\epsilon} F^2\,d\mu  &\le&  \frac{1}{\epsilon}\int_0^\epsilon\int_{\partial M} F^2(x,t)|\det D E(x,t)|d\umu dt\\
&\le&\frac{1}{\epsilon}\int_0^\epsilon\int_{\partial M} \left( |F(x,0)|+\int_0^t\left|\frac{\partial  F}{\partial s}(x,s)\right|ds\right)^2|\det DE(x,t)|d\umu dt\\ 
&\le&\frac{2}{\epsilon}\int_0^\epsilon\int_{\partial M}  F(x,0)^2|\det DE(x,t)|d\umu dt\\&&+\frac{2}{\epsilon}\int_0^\epsilon\int_{\partial M}\left(\int_0^t\left|\frac{\partial  F}{\partial s}(x,s)\right|ds\right)^2|\det DE(x,t)|d\umu dt\\
&\le&4\int_{\partial M}  F(x,0)^2d\umu +\frac{2}{\epsilon}\int_0^\epsilon\int_{\partial M}t\int_0^t\left|\frac{\partial  F}{\partial s}(x,s)\right|^2|\det DE(x,s)|\frac{|\det DE(x,t)|}{|\det DE(x,s)|}ds\,d\umu dt
\\&\le& 4\int_{\partial M} F^2d\umu +2\epsilon\int_{M_\epsilon} |\nabla  F|^2d\mu  \eq%
Taking $F=\un_{A_j}f$ in the above inequality  we get
\bqn{tracef}\frac{1}{\epsilon}\int_{A_{j,e}^\epsilon} f^2\,d\mu &\le& 4\int_{\partial_e A_{j}} {f}^2d\umu +2\epsilon\int_{A_{j}} |\nabla  {f}|^2d\mu\eqn

We proceed with bounding the second term $\epsilon \int_{A_j{\setminus A_{j,e}^\epsilon}} f^2\,d\mu $. 
Let $\xi:M\to \R_+$ be a Lipschitz function such that $|\nabla \xi|\le\frac{1}{\epsilon}$ and 
\[\begin{cases}\xi\equiv1\,,&\hbox{in $M{\setminus M^\epsilon}$}\\
0\le\xi\le1\,, &\hbox{in $M^\epsilon$}\\
\xi\equiv0\,,&\hbox{on $\partial M$}\end{cases}\] 
We get
\bq\epsilon \int_{A_j{\setminus A_{j,e}^\epsilon}} f^2\,d\mu  &\le&\epsilon \int_{A_j}\xi f^2\,d\mu=\epsilon \int_{M}\xi\un_{A_j} f^2\,d\mu\\
&{\le}& \epsilon P_1\int_{M} |\nabla (\xi \un_{A_j}f^2)|d\mu= \epsilon P_1\int_{A_j} |\nabla (\xi f^2)|d\mu\\  
&\le&\epsilon P_1\left(\int_{A_j}|\nabla\xi| f^2d\mu+2\int_{A_j}\xi f|\nabla f|d\mu\right)\\
&\le&\epsilon P_1\left(\frac{1}{\epsilon}\int_{A_{j,e}^\epsilon}f^2d\mu+2\left(\int_{A_j} (\xi f)^2d\mu\right)^{\frac{1}{2}}\left(\int_{A_j} |\nabla f|^2d\mu\right)^{\frac{1}{2}}\right) \\
&\stackrel{\eqref{tracef}}{\le}&\epsilon P_1\left(4\int_{\partial_e A_{j}} {f}^2d\umu +2\epsilon\int_{A_{j}} |\nabla  {f}|^2d\mu\right)\\
&&+2\epsilon P_1\bar\lambda_1(M)^{-1/2}\left(\int_{A_j} |\nabla (\xi f)|^2d\mu\right)^{\frac{1}{2}}\left(\int_{A_j}|\nabla f|^2d\mu\right)^{\frac{1}{2}}\\
&\le&4\epsilon P_1\int_{\partial_e A_{j}} {f}^2d\umu +2\epsilon^2 P_1\int_{A_{j}} |\nabla  {f}|^2d\mu\\
&&+2 P_1\bar\lambda_1(M)^{-1/2}\left(\sqrt{\epsilon}\left(\frac{1}{\epsilon}\int_{A_{j,e}^\epsilon}f^2d\mu\right)^{\frac 1 2}\left(\int_{A_j}|\nabla f|^2d\mu\right)^{\frac{1}{2}}+\epsilon\int_{A_j}|\nabla f|^2d\mu\right)\\
&\stackrel{\eqref{tracef}}{\le}&4\epsilon P_1\int_{\partial_e A_{j}} {f}^2d\umu +2\epsilon P_1(\epsilon+\bar\lambda_1(M)^{-{\frac1 2}})\int_{A_{j}} |\nabla  {f}|^2d\mu\\
&&+2\sqrt{\epsilon} P_1\lambda_1(M)^{-1/2}\left( 4\int_{\partial_e A_{j}} {f}^2d\umu +2\epsilon\int_{A_{j}} |\nabla  {f}|^2d\mu\right)^{\frac 1 2}\left(\int_{A_j}|\nabla f|^2d\mu\right)^{\frac{1}{2}}\\
&\le&4 \epsilon P_1\int_{\partial_e A_{j}} {f}^2d\umu+2\epsilon P_1 \left(\epsilon+(1+\sqrt{2})\bar\lambda_1(M)^{-{\frac1 2}}\right)\int_{A_{j}} |\nabla  {f}|^2d\mu\\
&&+4\sqrt{\epsilon} P_1\bar\lambda_1(M)^{-1/2}\left( \int_{\partial_e A_{j}} {f}^2d\umu \right)^{\frac 1 2}\left(\int_{A_j}|\nabla f|^2d\mu\right)^{\frac{1}{2}}\\
 \eq
 where $P_1$ is the $L^1$-Poincar\'e constant  and $\bar{\lambda}_1(M)$ is the first Dirichlet eigenvalue of  $M$.  In the second and fifth inequalities we used the Poincar\'e inequality on Sobolev spaces $W^{1,1}_0(M)$ and $W^{1,2}_0(M)$ respectively.  
Hence, for any $\epsilon\in(0,\epsilon_0)$ we get
\bq\frac{\int_{A_j} |\nabla f|^2d\mu  }{\int_{A_j}\rho_\epsilon f^2\,d\mu  }&\ge&\frac{\int_{A_j} |\nabla f|^2d\mu  }{4(1+\epsilon P_1)\int_{\partial_e A_j} f^2\,d\umu +C_1(\epsilon)\int_{A_{j}} |\nabla f|^2d\mu+C_2(\epsilon)\left( \int_{\partial_e A_{j}} {f}^2d\umu \right)^{\frac 1 2}\left(\int_{A_j}|\nabla f|^2d\mu\right)^{\frac{1}{2}} }\\&=&\psi_\epsilon\left(\frac{\int_{A_j} |\nabla f|^2d\mu  }{\int_{\partial_e A_j} f^2\,d\umu }\right)
\eq
where  $$C_1(\epsilon):=2\epsilon\left(1+ P_1 \left(\epsilon+(1+\sqrt{2})\bar\lambda_1(M)^{-{\frac1 2}}\right)\right),\qquad C_2(\epsilon):=4\sqrt{\epsilon} P_1\bar\lambda_1(M)^{-1/2}$$
and  $\psi_\epsilon:(0,\infty)\to(0,\infty)$ defined as
\bq\psi_\epsilon(u)&:=&\frac{u}{4(1+\epsilon P_1)+C_1(\epsilon)u+C_2(\epsilon)\sqrt{u}}\eq
is an increasing function. Remark that $\epsilon_0$ is independent of the set $A_j$ and depends only on $(M,g)$. Let $f_j$ be the eigenfunction associated with $\lambda_1^{\epsilon}(A_j)$.  
\bq\max_{j\in \llbracket k \rrbracket} \lambda_1^{\epsilon}(A_j)&=&\max_{j\in \llbracket k \rrbracket} \frac{\int_{A_j} |\nabla f_j|^2d\mu  }{\int_{A_j}\rho_\epsilon f_j^2\,d\mu  }\\&\ge&\max_{j\in \llbracket k \rrbracket}\psi_\epsilon\left(\frac{\int_{A_j} |\nabla f_j|^2d\mu  }{\int_{\partial_e A_j} f_j^2\,d\umu }\right)\\&\ge&\max_{j\in \llbracket k\rrbracket}\psi_\epsilon(\sigma_1(A_j))
=\psi_\epsilon(\max_{j\in \llbracket k\rrbracket}\sigma_1(A_j))\\&\ge&\psi_\epsilon(\inf_{(A_1,\cdots,A_k)\in\mathcal{A}_k}\max_{j\in \llbracket k\rrbracket}\sigma_1(A_j))\eq
Therefore,
{\bq\lim_{\epsilon\to0}\Lambda_k^\epsilon&\ge&\frac{1}{4} \kappa_k\eq}
which completes the proof.
\wwtbp\par
We  continue the proof of the theorem. By Theorem \ref{miclo}, we have 
\bq\lambda_k^\epsilon&\ge&\frac{c}{k^6}\Lambda_k^\epsilon\eq
Passing to the limit and applying Lemma \ref{keylem} and Theorem \ref{provenzano} we conclude:
\bq\sigma_k&=&\lim_{\epsilon\to0}\lambda_k^\epsilon\ge\frac{c}{k^6}\lim_{\epsilon\to0}\Lambda_k^\epsilon\ge\frac{c}{5k^6} \kappa_k\eq
\wwtbp\par
Similar to Propositions \ref{fs-1} and \ref{ms-1}, we have the following improvement on manifolds.\begin{pro} There is a universal positive constant $c'$  such that
\bq\fo k\in \N,\qquad \sigma_{2k}&\ge&
\frac{c'}{\log^2(k+1)}\iota_k\eq
\end{pro}
\proof
Due to \cite{MR2961569,MR3323580},  there is a universal positive constant $c_1$  such that
\bq\fo k\in \N,\qquad \lambda^\epsilon_{2k}&\ge&
\frac{c_1}{\log^2(k+1)}\Lambda^\epsilon_k\eq
Passing to the limite and using Lemmas  \ref{0stek} and \ref{keylem} we get  
\bq\fo k\in \N,\qquad \sigma_{2k}&\ge&
\frac{c_1}{4\log^2(k+1)}\kappa_k\ge\frac{c_1}{16\log^2(k+1)}\iota_k\eq
\wwtbp\par

\begin{rem} The methods and results above can be adapted to a more general Steklov eigenvalue problem 
\bq\left\{\begin{array}{ll}
{\rm div}( \phi\nabla f)=0\,,&\quad\hbox{in}\;\; M\\
\frac{\partial f}{\partial \nu}=\sigma \gamma f\,,&\quad\hbox{on} \;\; \partial M
\end{array}\right.\eq
where $\gamma$ is a continuous positive function on $\partial M$ and $\phi$ is a smooth positive function on $M$. 
But in this paper we stick to the so-called homogenous Steklov problem when $\phi=1$ and $\gamma=1$. 
\end{rem}

\begin{rem}We now give a more explicit  relationship between the higher order Cheeger constants and the higher order Cheeger--Steklov constants. Let $$\rho_k(M):=\inf_{(A_1,\cdots,A_k)\in\A_k}\max_{l\in\llbracket k\rrbracket}\rho(A_l)$$
We show that 
\bqn{ah-cheeger}\rho_k(M)=\inf_{(A_1,\cdots,A_k)\in\A_k}\max_{l\in\llbracket k\rrbracket}\eta(A_l)=:h_k(M)\eqn
where $h_k(M)$ denotes the $k$-th Cheeger constant. Indeed, it is easy to check that $ \rho(A)\le \eta(A)$ which implies  $\rho_k(M)\le h_k(M)$. Thus it is enough to show that for every $\epsilon>0$, we have $h_k(M)\le\rho_k(M)+\epsilon$. Note that  $$\fo B\subset A ,\qquad \rho(B)\ge\rho(A)$$ Recall the definition of $\rho(A)$ in \eqref{rhoa}. For every $\epsilon>0$, there exists $B\in \A$ subset of  $A$ such that $\bar B\cap \partial_i A=\emptyset$ and
\bqn{class}0\le\eta(B)-\rho(B)\le \eta(B)-\rho(A)<\epsilon\eqn   
 Let $\A_k^\epsilon$ be a subset  of $\A_k$ such that $$\fo  (A_1,\cdots, A_k)\in \A_k^\epsilon,\quad 0\le\eta(A_l)-\rho(A_l)<\epsilon,\quad \fo l\in\llbracket k\rrbracket$$
We claim    
$$\inf_{(A_1,\cdots,A_k)\in\A_k}\max_{l\in\llbracket k\rrbracket}\rho(A_l)=\inf_{(A_1,\cdots,A_k)\in\A_k^\epsilon}\max_{l\in\llbracket k\rrbracket}\rho(A_l)$$
Indeed, let \bq\left[(A_1,\cdots,A_k)\right]&:=&\left\{(\tilde A_1,\cdots,\tilde A_k)\in \A_k\,:\, \max_{l\in\llbracket k\rrbracket}\rho(A_l)=\max_{l\in\llbracket k\rrbracket}\rho(\tilde A_l)\right\}\eq
The definition of $\rho_k(M)$ does not change if we choose a representation in  each class $[(A_1,\cdots,A_k)]$ and   take infimum only over the family of  representations. By \eqref{class}, it is clear that each class has a representation in $\A_k^\epsilon$. This proves the claim. 
Therefore

$$\rho_k(M)=\inf_{(A_1,\cdots,A_k)\in\A_k^\epsilon}\max_{l\in\llbracket k\rrbracket}\rho(A_l)>\inf_{(A_1,\cdots,A_k)\in\A_k^\epsilon}\max_{l\in\llbracket k\rrbracket}\eta(A_l)-\epsilon
\ge h_k(M)-\epsilon$$
This proves identity \eqref{ah-cheeger}. 
Now for a given $(A_1,\cdots,A_k)\in\A_k$, let $l_{\max}\in\llbracket k\rrbracket$ be such that 
\bq\eta(A_{l_{\max}})&=&\max_{l\in\llbracket k\rrbracket}\eta(A_l)\eq 
Then we define 
\bq\bar h'_k(M)&\df&\inf_{(A_1,\cdots,A_k)\in\A_k}\rho'(A_{l_{\max}})\eq
It is easy to check that we have the  following lower bound for $\iota_k(M)$
\bqn{hcheeger}\iota_k(M)&\ge& h_k(M)\bar h'_k(M)\eqn

Similarly we can define  $$\rho'_k(M):=\inf_{(A_1,\cdots,A_k)\in\A_k}\max_{l\in\llbracket k\rrbracket}\rho'(A_l)$$
With the same argument as above,  the following equality holds. 
$$\rho'_k(M)=\inf_{(A_1,\cdots,A_k)\in\A_k}\max_{l\in\llbracket k\rrbracket}\eta'(A_l)=:h'_k(M)$$
For a given $(A_1,\cdots,A_k)\in\A_k$, let $l'_{\max}\in\llbracket k\rrbracket$ be such that $$\eta'(A_{l'_{\max}})=\max_{l\in\llbracket k\rrbracket}\eta'(A_l)$$ 
Then define
$$\bar h_k(M):=\inf_{(A_1,\cdots,A_k)\in\A_k}\rho(A_{l'_{\max}})$$
and we get
\[\iota_k(M)\ge \bar h_k(M) h'_k(M)\]
\end{rem}
Jammes in  \cite{MR3449183} considered  several examples  to show that for $k=2$ the geometric quantities  $\eta(B)$ and $\eta'(B)$ appearing in the definition of $\iota_k(M)$ are both necessary in the lower bound of $\sigma_2(M)$.  Inspired by his examples, we give examples which show the necessity of quantities such as  $\eta(B)$ and $\eta'(B)$ in the lower bound for all $k\in \N$.
 \begin{example}\label{ex1} Exemple 4 of \cite{MR3449183} can be  used to show the necessity of quantities such as $\eta(B)$ and $\eta'(B)$ in the definition of $\iota_k$ for all $k\ge 2$: Consider $M_m=N\times (0,L_m)$, where $N$ is a closed manifold and $L_m=\frac{1}{m}$. The Steklov spectrum of $M_m$ can be calculated explicitly, see \cite[Lemma 6.1]{MR2807105}. They are 
 \begin{equation}\label{ceg}\left\{\,0\,,{L_m}^{-1}, \sqrt{\lambda_k(N)}\tanh(\sqrt{\lambda_k(N)}L_m),\sqrt{\lambda_k(N)}\coth(\sqrt{\lambda_k(N)}L_m):\,k\in\N\,\right\}\end{equation}
 where $\lambda_k(N)$ are the Laplace eigenvalues of $N$.
   It is clear that for every $k\in\N$, $\sigma_k=O(\frac{1}{m})$ as $m\to\infty$, while $h_2(M_m)\ge c$ for some positive constant $c$ independent of $m$ as shown in \cite[Exemple 4]{MR3449183}.  Note that $h_k(M_m)$ is a non-decreasing sequence in $k$. Hence we have $h_k(M)\ge h_2(M_m)\ge c$, for every $k\ge2$. This together with \eqref{hcheeger} and Theorem \ref{cheeg-stek}  show the necessity of  a quantity such as $\eta'(B)$ in the definition of $\iota_k(M_m)$ for all $k\in\N$. 
   \end{example}
\begin{example}\label{ex2'}
Let $\SS^1$ be the unit circle and $\SS^1_m$ denote a circle of radius $m$ with their standard metric.  Consider the sequence $(M_m:={\SS}^1_m\times(0,m^{3/2}))_{m\in\NN}$ with product metric.  The set of Steklov eigenvalues $\sigma_k(M_m)$ is given by \eqref{ceg}. Note that $\lambda_k({\SS^1_m})=\frac{1}{m^2}\lambda_k({\SS^1})$. Hence, for any fixed $k\in \N$  we have $$\sigma_k(M_m)\sim m^{3/2}\lambda_k(\SS^1_m)=\frac{1}{\sqrt{m}}\lambda_k(\SS^1)\quad \text{as $m\to\infty$}$$ Therefore
$$ \fo  k\in\N,\qquad\lim_{m\to\infty}\sigma_k(M_m)=0$$ 
It is easy to check that for every $k\in\NN$, $\lim_{m\to\infty}h_k(M_m)=0$. Indeed, if we choose $A_l=S^1_m\times(\frac{(l-1)m^{3/2}}{k},\frac{l m^{3/2}}{k})$, $l\in\llbracket k\rrbracket$ then 
\bq h_k(M_m)&\le&\max_{l\in\llbracket k\rrbracket}\frac{\umu  (\partial_iA_l)}{\mu(A_l)}={\frac{4\pi m}{2\pi m^{5/2}/k}=\frac{2 k}{m^{3/2}}}\to0, \quad m\nearrow\infty\eq
We now show  that there exists a positive constant $C$ independent of $m$ such that $h'_k(M_m)\ge C$. Note that $h'_k(M_m)$  is a  non-decreasing sequence in $k$. Thus, it is enough to show that $h'_2(M_m)\ge C$ for some constant  $C>0$ independent of $m$.  Let $(A_1,A_2)$ be a partition of $M_m$ (w.l.o.g. we can assume $A_1$ is connected). Let assume $\partial_i A_1$ only intersect{s} one of the boundary components of $M_m$. Fixing the area of $A_1$,  $\max\left\{\frac{\umu_{m}  (\partial_iA_1)}{\umu_{m}  (\partial_eA_1)},\frac{\umu_{m}  (\partial_iA_2)}{\umu_{m}  (\partial_eA_2)}\right\}$ is minimized when $\partial _iA_1={\SS^1_m}\times\{x\}$ for some $x\in(0,m)$ {(where $\umu_m$ is the one-dimensional Riemannian measure of a set in $M_m$)}. Thus,
\bq1\le\max\left\{\frac{\umu_{m}  (\partial_iA_1)}{\umu_{m}  (\partial_eA_1)},\frac{\umu_{m}  (\partial_iA_2)}{\umu_{m}  (\partial_eA_2)}\right\}\eq
We now assume otherwise, i.e. $\partial_iA_1$ intersects both boundary components of $M_m$. We have\bq 
\max\left\{\frac{\umu_{m}  (\partial_iA_1)}{\umu_{m}  (\partial_eA_1)},\frac{\umu_{m}  (\partial_iA_2)}{\umu_{m}  (\partial_eA_2)}\right\}&\ge&\frac{2m^{\frac 3 2}}{2\pi m}=\frac{\sqrt{m}}{\pi}
\eq
We conclude that  for $m>\pi^2$,  $$h'_k(M_m)\ge h'_2(M_m)\ge1$$ This example shows the necessity of  a quantity such as $\eta(B)$ in the definition of $\iota_k(M_m)$ for all $k\in\N$.  For $k=2$, a similar example  has been studied in \cite{MR3449183}.    \\
\end{example}

\begin{example}[Cheeger dumbbell]
Girouard and   Polterovich in \cite{GP10} studied a family of Cheeger dumbbells  $M_\epsilon$ and showed   that $\lim_{\epsilon\to0} \sigma_k(M_\epsilon)=0$ for every $k\in\N$. In their example, $M_\epsilon$ is a  domain in $\R^2$ consists of the union of two Euclidean unit disks $\mathbb{D}_1\cup\mathbb{D}_2$ connected with a thin rectangular neck $L_\epsilon$ of length $\epsilon$ and width $\epsilon^3$. It is easy to check  that $h_2(M_\epsilon)\to 0$ as $\epsilon\to 0$. We show that 
   for $k\ge 3$,  $h_k(M_\epsilon)\ge c>0$, where   $c$  is a constant independent of $\epsilon$.   Since $h_k(M_\epsilon)\ge h_3(M_\epsilon)$, it is enough to show that $h_3(M_\epsilon)>c$.  By contrary, we assume that $h_3(M_\epsilon)\to0$ as $\epsilon\to 0$. Hence there is a family of $(A^\epsilon_1,A^\epsilon_2,A^\epsilon_3)$ such that $$\max\left\{\frac{\umu  (\partial_iA^\epsilon_1)}{\mu (A^\epsilon_1)},\frac{\umu  (\partial_iA^\epsilon_2)}{\mu (A^\epsilon_2)},\frac{\umu  (\partial_iA^\epsilon_3)}{\mu (A^\epsilon_3)}\right\}\to 0,\quad\epsilon\to0$$
   Hence we have $\partial_iA_l^\epsilon\subset L_\epsilon$,  for all $l\in\llbracket 3\rrbracket$. Therefore, there exists $l\in\llbracket 3\rrbracket$ such that $A_l^\epsilon\subset L_\epsilon$. (Notice that the argument   uses the fact that $M_\epsilon$ is a subset of  $\R^2$.) Taking $\epsilon=\frac{1}{m}$, $m\in \N$, and then using the similar argument as in  \cite[Exemple 4 ]{MR3449183}, we conclude that  for any $\epsilon$ small enough
\[\frac{\umu  (\partial_iA^\epsilon_l)}{\mu (A^\epsilon_l)}\ge c>0\]  
{where $c$ is independent of $\epsilon$}.
It is a contradiction. \\
This example as in Example \ref{ex1} shows the necessity of $\eta'(B)$ in $\iota_k(M_\epsilon)$. However,  in Example~\ref{ex1}  the volume of the family of manifolds tends to zero, while in this example the area and the boundary length of $M_\epsilon$ are uniformly controlled.

\end{example}

\appendix
\section{Continuity properties}\label{App}

In the context of Section \ref{MS}, it is natural to wonder if the mapping $P\mapsto K$  is continuous in some sense.
In particular it could provide an alternative approach to the acceleration technique 
in the deduction of Theorem \ref{theo1m} when $P$ has no  Dirichlet
gap on $M\setminus V$. We give here an example of a strong continuity result useful in this direction, but under quite restrictive assumptions of hyperboundedness on $K$ (holding e.g.\ when $\LL^2(\nu)$ is finite dimensional) and of ergodicity (see Lemma \ref{ergod} below).\par\me
For any $\epsilon\in(0,1)$, consider
\bq
P_\epsilon&\df& (1-\epsilon)P+\epsilon\mu\eq
where $\mu$ is seen as the Markov kernel given for any $x\in M$ by $\mu(x,dy)\df\mu(dy)$.
It is clear that $\mu$ is still reversible with respect to the Markov operator $P_\epsilon$. Its advantage is that $P_\epsilon$ has a Dirichlet gap on $M\setminus V$ larger than $\mu(V)\epsilon>0$ (consider the function $F\df \un_{M\setminus V}$ in its definition as an infimum).
Let $K_\epsilon$ be the Steklov Markov kernel associated to $P_\epsilon$ and $V$ (all notions relative to $P_\epsilon$ will receive $\epsilon$ in index).
A strong approximation property holds under certain circumstances:
\begin{lem}\label{appeps}
Let $\xi$ stands for the normalisation of the restriction of $\mu$ to $M\setminus V$ and let $\gamma$ be the law of $Z(\tau_V)$, where $(Z(n))_{n\in\ZZ_+}$
is a Markov chain starting with initial law $\xi$ and transition kernel $P$, and where $\tau_V$ is the hitting time of $V$.
Assume that \\
$\bullet$  The Radon-Nykodim density $d\gamma/d\nu$ belongs to $\LL^2(\nu)$.\\
$\bullet$ The operator $K$ is hyperbounded, in the sense there exist $1<q<2$ such that $K$ is bounded from $\LL^q(\nu)$ to $\LL^2(\nu)$.\par
Let $\vvv\cdot\vvv_\nu$ stands for the operator norm in $\LL^2(\nu)$, then we have
\bq
\lim_{\epsilon\ri 0_+} \vvv K_\epsilon-K\vvv_\nu&=&0\eq
\end{lem}
\proof
To simplify the notation with respect to \eqref{tauA}, introduce 
\bq
\tau&\df&\tau_{V}\\
\wi \tau&\df&\inf\{n\in\NN\St Z(n)\in V\}\eq
so that for any $f\in\LL^2(\nu)$, we have,
\bq
\hbox{$\nu$-a.s. in }x\in V,\qquad 
K[f](x)&=&\EE_x[f(Z(\wi \tau))]\eq
and similarly,
\bq
\hbox{$\nu$-a.s. in },x\in V\qquad 
K_\epsilon[f](x)&=&\EE_x[f(Z_{\epsilon}(\wi \tau_{\epsilon}))]\eq
The transitions of the Markov chain $Z_{\epsilon}$ can be interpreted in the following way: at each time,
the chain chooses with probability $1-\epsilon$ to go through a transition dictated by $P$ and with probability $\epsilon$
 a new position is sampled according to $\mu$.
It follows that for any $f\in\LL^2(\nu)$ and $x\in M$,
\bqn{tauwitau}
\EE_x[f(Z_{\epsilon}(\wi \tau_{\epsilon}))]&=&\EE_x[f(Z(\wi \tau))(1-\epsilon)^{\wi \tau}]+\EE_x[(1-(1-\epsilon)^{\wi \tau})]\int \EE_y[f(Z_{\epsilon}(\tau_{\epsilon}))]\,\mu(dy)\eqn
The last integral  is just $\EE_\mu[f(Z_{\epsilon}(\tau_{\epsilon}))]$ and
we remark that
\bq
\EE_\mu[f(Z_{\epsilon})(\tau_{\epsilon})]&=&\mu(V)\nu[f]+\mu(M\setminus V)\EE_\xi[f(Z_{\epsilon})(\wi \tau_{\epsilon})]\eq
To compute the last term, let us integrate \eqref{tauwitau} with respect to $\xi(dx)$, to get
\bq
\EE_\xi[f(Z_{\epsilon}(\wi \tau_{\epsilon}))]&=&\EE_\xi[f(Z(\wi \tau))(1-\epsilon)^{\wi \tau}]+\EE_\xi[(1-(1-\epsilon)^{\wi \tau})]\EE_\mu[f(Z_{\epsilon}(\tau_{\epsilon}))]\\
&=&\EE_\xi[f(Z(\wi \tau))(1-\epsilon)^{\wi \tau}]+\EE_\xi[(1-(1-\epsilon)^{\wi \tau})](\mu(V)\nu[f]+\mu(M\setminus V)\EE_\xi[f(Z_{\epsilon}(\wi \tau_{\epsilon}))])\eq
i.e.
\bqn{Cepf}
\EE_\xi[f(Z_{\epsilon})(\wi \tau_{\epsilon})]\ =\ C_\epsilon(f)\ \df\ \frac{\EE_\xi[f(Z(\wi \tau))(1-\epsilon)^{\wi \tau}]+\EE_\xi[(1-(1-\epsilon)^{\wi \tau})]\mu(V)\nu[f]}{1-\mu(M\setminus V)\EE_\xi[(1-(1-\epsilon)^{\wi \tau})]}\eqn
Putting together these computations, we get
\bq
K_\epsilon[f](x)&=&\EE_x[f(Z(\wi \tau))(1-\epsilon)^{\wi \tau}]+C_\epsilon(f)\EE_x[(1-(1-\epsilon)^{\wi \tau})]\eq
so that
\bq
K_\epsilon[f](x)-K[f](x)&=&\EE_x[f(Z(\wi \tau))((1-\epsilon)^{\wi \tau}-1)]+C_\epsilon(f)\EE_x[(1-(1-\epsilon)^{\wi \tau})]\eq
and
\bqn{vvv}
\nonumber\lVe K_\epsilon[f]-K[f]\rVe_{\LL^2(\nu)}^2&\leq &
2\int\EE^2_x[f(Z(\wi \tau))((1-\epsilon)^{\wi \tau}-1)]\,\nu(dx)\\&&+
2C_\epsilon(f)^2\int \EE^2_x[(1-(1-\epsilon)^{\wi \tau})]\,\nu(dx)\eqn
By using twice  H\" older's inequality with respect to the exponents $r,r'>1$ and $s,s'>1$ such that $1/r+1/r'=1$ and $1/s+1/s'=1$, we get that
\bq
\int\EE^2_x[f(Z(\wi \tau))((1-\epsilon)^{\wi \tau}-1)]\,\nu(dx)&\leq & \int\EE_x[f^r(Z(\wi \tau))]^{2/r}\EE_x[((1-\epsilon)^{\wi \tau}-1)^{r'}]^{2/r'}\,\nu(dx)\\
&=&\int (K[f^r](x))^{2/r}\EE_x[((1-\epsilon)^{\wi \tau}-1)^{r'}]^{2/r'}\,\nu(dx)\\
&\leq & 
\nu[(K[f^r])^{2s/r}]^{1/s}\lt(\int \EE_x[((1-\epsilon)^{\wi \tau}-1)^{r'}]^{2s'/r'}\,\nu(dx)\rt)^{1/s'}\eq
Due to the hyperboundedness assumption there exists $C\geq 1$ such that
\bq
\fo g\in\LL^{q}(\nu),\qquad \nu[(K[g])^2]^{1/2}&\leq & \nu[g^q]^{1/q}\eq
So with $g=f^r$, $r=2/q=s$, 
we get
\bq
\nu[(K[f^r])^{2s/r}]^{1/s}\ \leq \  (C\nu[f^2]^{1/q})^{2/r}\ =\ C^{2/r}\nu[f^2]\eq
and thus
with
\bq
r'\ =\ \frac{r}{r-1}\ =\ \frac{2}{2-q}\ =\ s'\eq
we deduce
\bq
\int\EE^2_x[f(Z(\wi \tau))((1-\epsilon)^{\wi \tau}-1)]\,\nu(dx)&\leq &
C^{2/r} \lt(\int \EE_x[((1-\epsilon)^{\wi \tau}-1)^{r'}]^{2}\,\nu(dx)\rt)^{1/r'}
\nu[f^2]\eq
By dominated convergence, we have 
\bq
\lim_{\epsilon\ri 0_+} \int \EE_x[((1-\epsilon)^{\wi \tau}-1)^{r'}]^{2}\,\nu(dx)&=&0\eq
thus we obtain
\bq
\lim_{\epsilon\ri 0_+}\sup_{f\in B(\LL^2(\nu))}\int\EE^2_x[f(Z(\wi \tau))((1-\epsilon)^{\wi \tau}-1)]\,\nu(dx)&=&0\eq
where $B(\LL^2(\nu))$ is the unitary ball of $\LL^2(\nu)$.\par
Since dominated convergence also implies
\bq
\lim_{\epsilon\ri 0_+}\int \EE^2_x[(1-(1-\epsilon)^{\wi \tau})]\,\nu(dx)&=&0\eq
to get the announced convergence, in view of \eqref{vvv}, it remains
to show that
\bq
\sup_{\epsilon\in(0,1),\, f\in B(\LL^2(\nu))} C_\epsilon^2(f)&<&+\iy\eq
Indeed, by the definition of $C_\epsilon(f)$ given in \eqref{Cepf}, we have
\bq
\fo \epsilon\in(0,1),\qquad
C_\epsilon^2(f)&\leq & 
2\frac{\EE^2_\xi[f(Z(\wi \tau))]+\mu^2(V)\nu^2[f]}{1-\mu(M\setminus V)}\\
&=&2\frac{\gamma^2[f]+\mu^2(V)\nu^2[f]}{\mu(V)}\\
&\leq &2\frac{\nu[(d\gamma/d\nu)^2]+\mu^2(V)}{\mu(V)}\nu[f^2]\eq
thus
\bq
\sup_{\epsilon\in(0,1),\, f\in B(\LL^2(\nu))} C_\epsilon^2(f)&\leq & 2\frac{\nu[(d\gamma/d\nu)^2]+\mu^2(V)}{\mu(V)}\eq
\wwtbp
\par
The first assumption of Lemma \ref{appeps} is linked to the square integrability 
of $\wi\tau\df\inf\{n\in\NN\St Z(n)\in V\}$ where $Z\df(Z(n))_{n\in\ZZ_+}$ is a Markov chain whose transition kernel is $P$
and initial distribution $\nu$:
\begin{lem}\label{ergod}
Assume the ergodic theorem holds for $Z$.
Then with the notation of Lemma \ref{appeps}, we have
\bq
\nu\lt[\lt(\frac{d\gamma}{d\nu}\rt)^2\rt]&\leq & \frac{\sqrt{\EE_\nu[(\wi\tau-1)^2]}}{\EE_\nu[\wi\tau-1]}\eq
with the convention that the r.h.s.\ is $+\iy$ when $\EE_\nu[\wi\tau]=+\iy$.
\end{lem}
\proof
Of course, it is sufficient to consider the case where $\EE_\nu[\wi\tau]<+\iy$.
Then for any measurable and bounded function $g$ on $M\setminus V$, we have
\bqn{ergo}
\EE_\nu\lt[\sum_{n=1}^{\wi\tau-1} g(Z(n))\rt]&=&\mu(M\setminus V)\EE_\nu[\wi\tau]\xi[g]\eqn
This is a classical consequence of the ergodic theorem: consider the sequence of passage times by $V$ defined by the iteration
\bq
\wi\tau_0&\df&0\\
\fo p\in\ZZ_+,\qquad \wi\tau_{p+1}&\df& \inf\{n> \wi\tau_p\St Z(n)\in V\}\eq
and the process $Y\df (Y(p))_{p\in\ZZ_+}$ trace of $Z$ on $V$:
\bq
\fo p\in\ZZ_+,\qquad Y(p)&\df& Z(\wi\tau_p)\eq
Note that $Y$ is a stationary Markov chain whose transition are given by $K$, leaving $\nu$ invariant.
Consider the mapping $G$ defined on $V$ by
\bq
\fo x\in V,\qquad 
G(x)&\df& \EE_x\lt[\sum_{n=1}^{\wi\tau-1} g(Z(n))\rt]\eq
and extend $g$ on $M$ by making it vanish on $V$.
Then we have a.s.
\bq
\lim_{p\ri\iy} \frac1{p+1}\sum_{n=0}^{p}G(Y_n)&=&\nu[G]\\
\lim_{p\ri\iy} \frac{1}{\wi\tau_{p+1}}\sum_{n=0}^{\wi\tau_{p+1}-1}g(Z_n)&=&
\mu[g]\ =\ \mu(M\setminus V)\xi[g]\\
\lim_{p\ri\iy} \frac{\wi\tau_{p+1}}{p+1}&=&\EE_\nu[\wi\tau]\eq
The relation \eqref{ergo} then follows from the equality
\bq
\fo p\in\ZZ_+,\qquad \sum_{n=0}^{p}G(Y_n)&=&\sum_{n=0}^{\wi\tau_{p+1}-1}g(Z_n)\eq
\par
Taking $g=\un_{M\setminus V}$, we get $\EE_\nu[\tau-1]=\mu(M\setminus V)\EE_\nu[\tau]$  in  \eqref{ergo}, which can be written under the form
\bq
\EE_\nu\lt[\sum_{n=1}^{\wi\tau-1} g(Z(n))\rt]&=&\EE_\nu[\wi\tau-1]\xi[g]\eq
\par
Consider $f$ a measurable and bounded function on $V$ and associate the mapping $g$ defined on $M\setminus V$ by
\bq
\fo x\in M\setminus V,\qquad 
g(x)&\df& \EE_x[f(Z(\wi\tau))]\eq
 On one hand, we have by definition,
\bq
\gamma[f]\ =\ \EE_\xi[f(Z(\wi\tau))]\ =\ \xi[g]\ =\ \frac{\EE_\nu\lt[\sum_{n=1}^{\wi\tau-1} g(Z(n))\rt]}{\EE_\nu[\wi\tau-1]}\eq
On the other hand, we compute that
\bq
\EE_\nu\lt[\sum_{n=1}^{\wi\tau-1} g(Z(n))\rt]&=&\sum_{n\in\NN}\EE_\nu\lt[\un_{n<\wi\tau}  g(Z(n))\rt]\\
&=&\sum_{n\in\NN}\EE_\nu\lt[\un_{n<\wi\tau}  \EE_{Z(n)}\lt[f(Z(\wi\tau))\rt]\rt]\\
&=&\sum_{n\in\NN}\EE_\nu\lt[\un_{n<\wi\tau}  f(Z(\wi\tau))\rt]\eq
where we used the Markov property and the fact that $\{n<\wi\tau\}$ belongs to the $\sigma$-field generated by $(Z(p))_{p\in\lin n\rin}$.
The last sum is equal to
\bq
\EE_\nu\lt[\sum_{n\in\NN}\un_{n<\wi\tau}  f(Z(\wi\tau))\rt]&=&\EE_\xi\lt[(\wi\tau-1) f(Z(\wi\tau))\rt]\eq
By Cauchy-Schwarz' inequality the last term is bounded
by 
\bq
\lve \EE_\nu\lt[(\wi\tau-1) f(Z(\wi\tau))\rt]\rve &\leq & \sqrt{\EE_\nu[(\wi\tau-1)^2]\nu[(K[f])^2]}\\
&\leq &\sqrt{\EE_\nu[(\wi\tau-1)^2]\nu[f^2]} \eq
Putting these observations together, we get
\bq
\lve \gamma[f]\rve &\leq & \frac{\sqrt{\EE_\nu[(\wi\tau-1)^2]}}{\EE_\nu[\wi\tau-1]}\sqrt{\nu[f^2]}\eq
first for any $f$  measurable and bounded function on $V$ and next by completion for all $f\in\LL^2(\nu)$.
The announced result follows by the Hilbert space's duality.
\wwtbp


\vskip2cm
\hskip70mm
\vbox{
\copy5
\vskip5mm
\copy6
}
\end{document}